  \documentclass[12pt,a4paper]{article}

  \makeatletter
  \@addtoreset{equation}{section}
  \makeatother

    \makeatletter
    \@addtoreset{figure}{section}
    \makeatother

    \makeatletter
    \@addtoreset{table}{section}
    \makeatother

\usepackage[]{hyperref}

\usepackage{color}

\definecolor{SteelBlue}  {rgb}{0.275,0.51,0.706}

\usepackage{latexsym}
\usepackage{ifthen}

\usepackage{amsmath}
\usepackage{amsthm}
\usepackage{amscd}
\usepackage{amsfonts}
\usepackage{amssymb}

\usepackage{ucs}
\usepackage[utf8]{inputenc}
\usepackage{longtable} 
\usepackage{booktabs} 
\usepackage[normalem]{ulem}
\usepackage{tabularx} 
\usepackage{float}

\usepackage{colortbl}

 \newtheoremstyle{TStyle_acknowledgement}{1.25em}{1.25em}{\itshape}{}{\bfseries\upshape}{.}{ }{}
 \theoremstyle{TStyle_acknowledgement}
 
 \newtheoremstyle{TStyle_algorithm}{1.25em}{1.25em}{\itshape}{}{\bfseries\upshape}{.}{ }{}
 \theoremstyle{TStyle_algorithm}
 \newtheorem{algorithm}{Algorithm}[section]
 \newtheoremstyle{TStyle_axiom}{1.25em}{1.25em}{\itshape}{}{\bfseries\upshape}{.}{ }{}
 \theoremstyle{TStyle_axiom}
 
 \newtheoremstyle{TStyle_case}{1.25em}{1.25em}{\itshape}{}{\bfseries\upshape}{.}{ }{}
 \theoremstyle{TStyle_case}
 
 \newtheoremstyle{TStyle_claim}{1.25em}{1.25em}{\itshape}{}{\bfseries\upshape}{.}{ }{}
 \theoremstyle{TStyle_claim}
 
 \newtheoremstyle{TStyle_comment}{1.25em}{1.25em}{\itshape}{}{\bfseries\upshape}{.}{ }{}
 \theoremstyle{TStyle_comment}
 
 \newtheoremstyle{TStyle_conclusion}{1.25em}{1.25em}{\itshape}{}{\bfseries\upshape}{.}{ }{}
 \theoremstyle{TStyle_conclusion}
 
 \newtheoremstyle{TStyle_condition}{1.25em}{1.25em}{\itshape}{}{\bfseries\upshape}{.}{ }{}
 \theoremstyle{TStyle_condition}
 
 \newtheoremstyle{TStyle_conjecture}{1.25em}{1.25em}{\itshape}{}{\bfseries\upshape}{.}{ }{}
 \theoremstyle{TStyle_conjecture}
 
 \newtheoremstyle{TStyle_corollary}{1.25em}{1.25em}{\itshape}{}{\bfseries\upshape}{.}{ }{}
 \theoremstyle{TStyle_corollary}
 
 \newtheoremstyle{TStyle_criterion}{1.25em}{1.25em}{\itshape}{}{\bfseries\upshape}{.}{ }{}
 \theoremstyle{TStyle_criterion}
 
 \newtheoremstyle{TStyle_definition}{1.25em}{1.25em}{\itshape}{}{\bfseries\upshape}{.}{ }{}
 \theoremstyle{TStyle_definition}
 
 \newtheoremstyle{TStyle_example}{1.25em}{1.25em}{\itshape}{}{\bfseries\upshape}{.}{ }{}
 \theoremstyle{TStyle_example}
 
 \newtheoremstyle{TStyle_exercise}{1.25em}{1.25em}{\itshape}{}{\bfseries\upshape}{.}{ }{}
 \theoremstyle{TStyle_exercise}
 
 \newtheoremstyle{TStyle_lemma}{1.25em}{1.25em}{\itshape}{}{\bfseries\upshape}{.}{ }{}
 \theoremstyle{TStyle_lemma}
 \newtheorem{lemma}[algorithm]{Lemma}
 \newtheoremstyle{TStyle_notation}{1.25em}{1.25em}{\itshape}{}{\bfseries\upshape}{.}{ }{}
 \theoremstyle{TStyle_notation}
 
 \newtheoremstyle{TStyle_problem}{1.25em}{1.25em}{\itshape}{}{\bfseries\upshape}{.}{ }{}
 \theoremstyle{TStyle_problem}

 \newtheoremstyle{TStyle_proposition}{1.25em}{1.25em}{\itshape}{}{\bfseries\upshape}{.}{ }{}
 \theoremstyle{TStyle_proposition}
 \newtheorem{proposition}[algorithm]{Proposition}
 \newtheoremstyle{TStyle_remark}{1.25em}{1.25em}{\itshape}{}{\bfseries\upshape}{.}{ }{}
 \theoremstyle{TStyle_remark}
 
 \newtheoremstyle{TStyle_solution}{1.25em}{1.25em}{\itshape}{}{\bfseries\upshape}{.}{ }{}
 \theoremstyle{TStyle_solution}
 
 \newtheoremstyle{TStyle_summary}{1.25em}{1.25em}{\itshape}{}{\bfseries\upshape}{.}{ }{}
 \theoremstyle{TStyle_summary}
 
 \newtheoremstyle{TStyle_theorem}{1.25em}{1.25em}{\itshape}{}{\bfseries\upshape}{.}{ }{}
 \theoremstyle{TStyle_theorem}
 \newtheorem{theorem}[algorithm]{Theorem}

\usepackage{sectsty}

   \allsectionsfont{\setlength{\parindent}{0em}\textcolor{SteelBlue}}

\ifx\pdfoutput\undefined 
  \usepackage[dvips]{graphicx}
\else 
  \usepackage{graphicx}
\fi
\newcommand{\includegraphicsX}[2]{
\newdimen\imgwidth
\setbox0=\hbox{\includegraphics[#1]{#2}}
\imgwidth=\wd0
\ifthenelse{\lengthtest{\imgwidth > \textwidth}}
  {
     {\includegraphics[width=\textwidth]{#2}}
  }
  {
     {\includegraphics[#1]{#2}}
  }
}

\newcommand{\tableX}[3]{
\newdimen\tablewidth
\setbox1=\hbox{\begin{tabular}{#1}#3\end{tabular}}
\tablewidth=\wd1
\ifthenelse{\lengthtest{\tablewidth > \textwidth}}
  {
     \begin{tabularx}{\textwidth}{#2}
     #3
     \end{tabularx}
  }
  {
     \begin{tabular}{#1}
     #3
     \end{tabular}
  }
}

\newcommand{\EmptyParagraph}{ $\phantom{X}$}

\makeatletter
\let\new@ifnextchar\@ifnextchar
\let\new@ifnch\@ifnch
\makeatother

\def\thebibliography#1{
 \list
 {[\arabic{enumi}]}{\settowidth\labelwidth{[#1]}\leftmargin\labelwidth
 \advance\leftmargin\labelsep
 \usecounter{enumi}}
 \def\newblock{\hskip .11em plus .33em minus -.07em}
 \sloppy\clubpenalty4000\widowpenalty4000
 \sfcode`\.=1000\relax}

\makeatletter
\newcommand{\xleftrightarrows}[2][]{\mathrel{%
 \raise.40ex\hbox{$\ext@arrow 3095\leftarrowfill@{\phantom{#1}}{#2}$}%
 \setbox0=\hbox{$\ext@arrow 0359\rightarrowfill@{#1}{\phantom{#2}}$}%
 \kern-\wd0 \lower.40ex\box0}}

\newcommand{\xrightleftarrows}[2][]{\mathrel{%
 \raise.40ex\hbox{$\ext@arrow 3095\rightarrowfill@{\phantom{#1}}{#2}$}%
 \setbox0=\hbox{$\ext@arrow 0359\leftarrowfill@{#1}{\phantom{#2}}$}%
 \kern-\wd0 \lower.40ex\box0}}
\newcommand{\xleftrightarrow}[2][]{%
     \ext@arrow 0055{\leftrightarrowfill@}{#1}{#2}%
}
\def\leftrightarrowfill@{%
 \arrowfill@\leftarrow\relbar\rightarrow%
}
\makeatother

\begin{document}

\begin{center}
        \noindent {\large A new family of periodic functions as explicit roots of }
\end{center}
     
\begin{center}
        \noindent {\large a class of polynomial equations}
\end{center}
     
                 \begin{table}[H]
                 \begin{center}
                  
          {\renewcommand{\arraystretch}{1.1}
    \tableX{l}{X}{\multicolumn{1}{c}{Marc Artzrouni} \\ \multicolumn{1}{c}{Department of Mathematics} \\ \multicolumn{1}{c}{University of¨ Pau} \\ \multicolumn{1}{c}{64000 Pau} \\ \multicolumn{1}{c}{FRANCE} \\ }}  
        \end{center}
      \end{table}
      
\begin{center}
        \noindent \emph{Marc.Artzrouni@univ-pau.fr}
\end{center}
     
\begin{center}
        \noindent \emph{October 22, 2004}
\end{center}

       \EmptyParagraph
     
    \begin{abstract}
     
        \noindent For any positive integer $ n , \hspace{2mm} $a new family of periodic functions in power series form and of period $ n $ is used to solve in closed form a class of polynomial equations of order $ n.\hspace{2mm} $The $ n\hspace{2mm} $roots are the values of the appropriate function from that family taken at $ 0 , 1 , ... , n - 1. $
     
       \end{abstract}

      \section{Background}\label{SECTION.d8ca7572-f45a-479a-be71-9040b5475d5a}

        \noindent For centuries mathematicians have sought closed-form expressions for the roots of polynomial equations of arbitrary order $ n.\hspace{2mm} $It is now well known that polynomials of order four or less can be solved explicitly using rational operations and finite root extractions. The Abel-Ruffini theorem shows that this cannot be done for orders five and above  \cite{BIBITEM.32b8409e-4efa-45ca-ae05-adb468ace982}.

       Several authors have proposed series solutions of algebraic (and polynomial) equations ( \cite{BIBITEM.8d478637-3e3c-4e3c-b709-d4353976b5ee},  \cite{BIBITEM.09654ad0-cb4d-4f5e-97a6-455e481aaa20},  \cite{BIBITEM.1e4a0b01-ced4-402b-a8fa-48c735e80cb1},  \cite{BIBITEM.2ce5749f-bb1c-4c48-9c60-73b87c19a339},  \cite{BIBITEM.54e711a8-276d-4617-bd00-cf34c8750ced}). These solutions, which rely on hypergeometric functions, are cumbersome to implement and have not provided feasible alternatives to standard numerical methods.

       For a positive integer $ n , \hspace{2mm} $we will describe here in elementary fashion a family of functions

       \begin{equation}\label{EQUATION.e4267f21-c149-484a-995b-9d3574d58deb}
       x{ \left( t , u \right) } = { \sum }_{{m = 1}}^{{ \infty }}{ \beta }_{{m}}e^{{2u \pi m \times i / n}} \times t^{{m}}
       \end{equation}

        \noindent of two real variables $ t $ and $ u\hspace{2mm} $that are power series in $ t $ and are periodic of period $ n $ in the variable $ u.\hspace{2mm} $(These functions are analogous to complex Fourier series except that the summation extends over positive indices only).

       For a limited class of polynomial equations of order $ n\hspace{2mm} $(parameterized in some way by $  \left. t \right) \hspace{2mm} $the $ { \beta }_{{m}} ' s\hspace{2mm} $ will be found explicitly so that the series $  \sum { \left| { \beta }_{{m}}t^{{m}} | \hspace{2mm} \right. } $converges and the values of $ x{ \left( t , u \right) }\hspace{2mm} $at $ t\hspace{2mm} $and $ u = 0 , 1 , ... , n - 1 $ will be the $ n $ roots of the equation. The function $ x{ \left( t , u \right)  , }\hspace{2mm} $which could be thought as "elementary" in the same way

       \begin{equation}\label{}
       cos{ \left( x \right) } = {{ \sum }_{{n = 0}}^{{ \infty }}}{{ \left(  - 1 \right) }}^{{n}}x^{{2n}} / { \left( 2n \right) } ! \hspace{2mm}
       \end{equation}

        \noindent is an elementary function, thus provides explicit solutions simply through its values at $ t\hspace{2mm} $and $ u =  $$ 0 , 1 , ... , n - 1 $.

       This will be only a first step as the class of polynomial equations solved explicitly with these power series is limited. The ultimate goal is to generalize the approach proposed here to any polynomial equation.

       \EmptyParagraph

      \section{Preliminaries}\label{SECTION.df7b0848-27e7-4e7b-9a00-1cdecc717719}

        \noindent We start off with a polynomial equation in the form

       \begin{equation}\label{EQUATION.07dd3654-31ac-4ca3-a570-c25959e51126}
       x^{{n}} = {{a}_{{{n - 1}}}x^{{n - 1}} + a_{{n - 2}}x^{{n - 2}} + ... + a_{{1}}x + {a}_{{0}}}
       \end{equation}

        \noindent where the $ a_{{k}} ' s $ are real or complex coefficients and $ {a}_{{0}} =  \rho {e}^{{i \theta }}\hspace{2mm} $($  \left.  \rho  > 0 ,  -  \pi  <  \theta  \le  \pi  \right) \hspace{2mm} $is assumed throughout to be non-zero. (Otherwise  \eqref{EQUATION.07dd3654-31ac-4ca3-a570-c25959e51126} can trivially be reduced to an equation of degree $  \left. n - 1 \right) .\hspace{2mm}\hspace{2mm} $

       Equation  \eqref{EQUATION.07dd3654-31ac-4ca3-a570-c25959e51126} is transformed by multiplying the right-hand side by $ t^{{n\hspace{2mm}}} $where $ t $ is a real variable that we may initially think of as small but is destined to take on any real value including 1:

       \EmptyParagraph

       \begin{equation}\label{EQUATION.138330db-ce67-4d93-a8f0-635652033f0e}
       x^{{n}} = { \left( {a}_{{{n - 1}}}x^{{n - 1}} + a_{{n - 2}}x^{{n - 2}} + ... + a_{{1}}x + {a}_{{0}} \right) }t^{{n}}.
       \end{equation}

       \EmptyParagraph

       We will find a solution in the form of a function $ x{ \left( t , u \right) } $ given in Eq.  \eqref{EQUATION.e4267f21-c149-484a-995b-9d3574d58deb}. We will show that $  \sum { \left| { \beta }_{{m}}t^{{m}} \right| } $ converges (and $ x{ \left( t , u \right) } $ is therefore defined) when $ t\hspace{2mm} $is small enough or $ { \left| a_{{0}} | \hspace{2mm} \right. } $large enough, or the modulii $ { \left| a_{{k}} | \hspace{2mm} \right. } $other than $ { \left| a_{{0}} | \hspace{2mm} \right. } $small enough.

       When the series converges the values of $ x{ \left( t , u \right) } $ at $ t $ and $ u = 0 , 1 , ... , n - 1 $ will provide the $ n $ roots of Equation  \eqref{EQUATION.138330db-ce67-4d93-a8f0-635652033f0e}. Indeed, we will show that if we define the partial sums

       \begin{equation}\label{EQUATION.6e641714-973c-4278-a36b-890fb646588e}
       x{{ \left( t , u \right) }}_{{q}}\overset{{def.}}{{ = }}\underset{{m = 1}}{\overset{{q}}{ \sum }}{ \beta }_{{m}}e^{{2u \pi m \times i / n}}{{ \times t}}^{{m}} , \hspace{2mm}
       \end{equation}

        \noindent then for $ u = 0 , 1 , ... , n - 1\hspace{2mm} $the polyonomial in $ t $

       \begin{equation}\label{EQUATION.d182e62d-c715-4015-8a4a-f408167c8ca0}
       E{ \left( x{{ \left( t , u \right) }}_{{q}} \right) }\overset{{def.}}{{ = }}{{{ \left[ x{{ \left( t , u \right) }}_{{q}} \right] }}^{{n}} - { \left( {a}_{{{n - 1}}}{{ \left[ x{{ \left( t , u \right) }}_{{q}} \right] }}^{{n - 1}} + a_{{n - 2}}{{ \left[ x{{ \left( t , u \right) }}_{{q}} \right] }}^{{n - 2}} + ... + a_{{1}}x{{ \left( t , u \right) }}_{{q}} + {a}_{{0}} \right) }t^{{n}}}
       \end{equation}

        \noindent will approach $ 0 $ for $ q \rightarrow  \infty . $

       We begin by seeking the roots expressed as the infinite series

       \begin{equation}\label{EQUATION.c5df2467-e55b-4766-9170-787c0bdffe48}
       x{{ \left( t \right) } = b_{{1}}t + b_{{2}}t^{{2}} + b_{{3}}t^{{3}} + ...}
       \end{equation}

        \noindent where \emph{a priori} we will need $ n $ different sequences $ { \left\{ b_{{m}} \right\} }\hspace{2mm} $to generate the $ n $ roots. We will find these $ b_{{m}} ' s\hspace{2mm} $and show that in fact they are each of the form

       \begin{equation}\label{EQUATION.2837254a-062d-4249-928b-39562f5dadd4}
       b_{{m}} = { \beta }_{{m}}e^{{2k \pi m \times i / n}} , k = 0 , 1 , ... , n - 1
       \end{equation}

        \noindent required in Eq.  \eqref{EQUATION.e4267f21-c149-484a-995b-9d3574d58deb} with a single sequence $ { \left\{ { \beta }_{{m}} \right\} }\hspace{2mm} $that will be obtained explicitly through a simple discrete dynamical system.

       \EmptyParagraph

       The equation to solve is now

       \begin{equation}\label{EQUATION.4e36cc79-02db-45e8-86b0-cff2382aeaa5}
       {{x{ \left( t \right) }}}^{{n}} = { \left( {a}_{{{n - 1}}}{{x{ \left( t \right) }}}^{{n - 1}} + a_{{n - 2}}{{x{ \left( t \right) }}}^{{n - 2}} + ... + a_{{1}}x{ \left( t \right) } + {a}_{{0}} \right) }t^{{n}}.
       \end{equation}

        \noindent Before proceeding we need some notations and preliminary results.

       \EmptyParagraph

       We define the powers $ {B}_{{q}}^{{r}} $ of the partial sums of $ x{ \left( t \right) } $:

       \begin{equation}\label{}
       {B}_{{q}}^{{r}}\overset{{def.}}{{ = }}{{ \left( b_{{1}}t + b_{{2}}t^{{2}} + {{... + b}}_{{q}}t^{{q}} \right) }}^{{r}},\hspace{2mm}q,r \in \mathrm{ \mathbb{N} } , 
       \end{equation}

        \noindent where for ease of notation the functional dependence of $ {B}_{{q}}^{{r}} $ on $ t $ is omitted. We let $ K{ \left( d,{B}_{{q}}^{{r}} \right) } $ denote the coefficient of $ t^{{d}}\hspace{2mm} $in $ {B}_{{q}}^{{r}} $ . ($ K{ \left( d,{B}_{{q}}^{{r}} \right) } $ does not depend on $ t $ and $  \left. K{ \left( d,{B}_{{q}}^{{r}} \right) } = 0\hspace{2mm}if\hspace{2mm}d < r\hspace{2mm}or\hspace{2mm}d > r.q \right) . $
     
       \begin{proposition}\label{PROPOSITION.118152ec-f39f-4537-96a3-9f97cc509b6a}
     
        \noindent With $ d \ge r , \hspace{2mm} $the coefficients $ K{ \left( d,{B}_{{q}}^{{r}} \right) } $ satisfy$ \hspace{2mm} $

       \begin{equation}\label{EQUATION.0269eb99-aa8d-4ad2-921d-7baf70e590b2}
       {K{ \left( d,{B}_{{q}}^{{r}} \right)  = K{ \left( d,{B}_{{d - r + 1}}^{{r}} \right) }\hspace{2mm} \forall q \ge d - r + 1,}}
       \end{equation}

        \noindent and

       \begin{equation}\label{EQUATION.6235dca5-7ffb-4919-8483-b028ecf45f0a}
       K{ \left( d,{B}_{{q}}^{{r}} \right) }{ = \underset{{m = 1}}{\overset{{min{ \left( d - r + 1 , q \right) }}}{ \sum }}{b}_{{m}}K{ \left( d - m,{ \left. {B}_{{q}}^{{r - 1}} \right) },\hspace{2mm}\hspace{2mm}r \ge 2 , \hspace{2mm}d \ge r ,  \right. }}
       \end{equation}

       \begin{equation}\label{EQUATION.453f5f3b-86ea-4183-9ce2-6adf54a01a37}
       K{ \left( q + s,{B}_{{q + s}}^{{s + 1}} \right) }{ = b_{{1}}K{ \left( q + s - 1,{B}_{{q + s - 1}}^{{s}} \right) } + \underset{{m = 2}}{\overset{{q - 1}}{ \sum }}{b}_{{m}}K{ \left( q + s - m,{B}_{{q + s - m}}^{{s}} \right)  + {b}_{{q}}{b}_{{1}}^{{s}}.}}
       \end{equation}

       \EmptyParagraph
     
    \end{proposition}\begin{proof}
     
        \noindent Equation  \eqref{EQUATION.0269eb99-aa8d-4ad2-921d-7baf70e590b2} expresses the fact that when $ q \ge  $ $ d - r + 1\hspace{2mm} $then only the first $ d - r + 1\hspace{2mm} $$ {b}_{{i}} ' s $ enter into $ K{ \left( d,{B}_{{q}}^{{r}} \right) }.\hspace{2mm} $Equation  \eqref{EQUATION.6235dca5-7ffb-4919-8483-b028ecf45f0a} is the convolution rule used to express the coefficient of $ {t}^{{d}}\hspace{2mm}\hspace{2mm} $in $ {B}_{{q}}^{{r}} $ considered as the product $ {B}_{{q}}^{{r - 1}} \times {B}_{{q}}^{{1}}. $ Equation  \eqref{EQUATION.453f5f3b-86ea-4183-9ce2-6adf54a01a37} uses  \eqref{EQUATION.0269eb99-aa8d-4ad2-921d-7baf70e590b2} to express a form of Eq.  \eqref{EQUATION.6235dca5-7ffb-4919-8483-b028ecf45f0a} in which all the $ K{ \left( a,{B}_{{u}}^{{v}} \right)  ' s} $ have identical values for $ a $ and $ u $. 
     
       \end{proof}
    
       \EmptyParagraph

       With $ x{ \left( t \right) } $ given in Eq.  \eqref{EQUATION.c5df2467-e55b-4766-9170-787c0bdffe48} both sides of Eq.  \eqref{EQUATION.4e36cc79-02db-45e8-86b0-cff2382aeaa5} are polynomials in $ t\hspace{2mm} $with powers $  \ge n. $ The goal is to find recursively the $ b_{{m}} ' s\hspace{2mm} $so that for each $ p \ge n $ the coefficients of $ t^{{p}}\hspace{2mm} $are equal on both sides of  \eqref{EQUATION.4e36cc79-02db-45e8-86b0-cff2382aeaa5}. This will be done by considering the roots $ x{ \left( t \right) } $ in the form of the gradually expanding partial sums $ {B}_{{q}}^{{1}}.\hspace{2mm} $ The coefficient $ b_{{1}}\hspace{2mm} $will be found by seeking a solution of the form $ {B}_{{1}}^{{1}}\hspace{2mm} $for which the coefficients of $ t^{{n}}\hspace{2mm} $on both sides of  \eqref{EQUATION.4e36cc79-02db-45e8-86b0-cff2382aeaa5} are equal. With $ b_{{1}}\hspace{2mm} $thus determined, $ b_{{2}}\hspace{2mm} $is found by seeking a solution of the form $ {B}_{{2}}^{{1}}\hspace{2mm} $for which the coefficients of $ t^{{n + 1}}\hspace{2mm} $on both sides of  \eqref{EQUATION.4e36cc79-02db-45e8-86b0-cff2382aeaa5} are equal, etc.

       We begin with a candidate solution $ {B}_{{1}}^{{1}} = b_{{1}}t $ by setting equal the coefficients of $ t^{{n}} $ on both side of Eq.  \eqref{EQUATION.4e36cc79-02db-45e8-86b0-cff2382aeaa5}, i.e.

       \begin{equation}\label{}
       B{ \left( n , {B}_{{1}}^{{n}} \right) } = {b}_{{1}}^{{n}}{{ = a}}_{{0}}.
       \end{equation}

        \noindent Therefore the $ n\hspace{2mm} $possible values of $ b_{{1}}\hspace{2mm} $are

       \begin{equation}\label{EQUATION.455ad3f5-5524-449d-91f1-acde5a424559}
       b_{{1}} = {{ \rho }}^{{1 / n}}{e}^{{i{ \left(  \theta  + 2k \pi  \right)  / n}}}\hspace{2mm}k = 0 , 1 , ... , n - 1\hspace{2mm}
       \end{equation}

        \noindent which is Eq.  \eqref{EQUATION.2837254a-062d-4249-928b-39562f5dadd4} with $ m = 1\hspace{2mm} $and

       \begin{equation}\label{}
       { \beta }_{{1}}\overset{{def.}}{{ = }}{{ \rho }}^{{1 / n}}{e}^{{i{ \theta  / n}}}.\hspace{2mm}
       \end{equation}

        \noindent We note that the $ n $ values $ x{ \left( t \right)  = b_{{1}}t\hspace{2mm}} $ with $ b_{{1}}\hspace{2mm} $of Eq.  \eqref{EQUATION.455ad3f5-5524-449d-91f1-acde5a424559} are the exact trivial solutions of Eq. \eqref{EQUATION.4e36cc79-02db-45e8-86b0-cff2382aeaa5} when all coefficients other than $ a_{{0}} $ are 0; otherwise they provide crude first-approximation solutions when $ t $ is small.

       $ \hspace{2mm} $

       In order to find $ b_{{2}}\hspace{2mm} $with a candidate solution $ {B}_{{2}}^{{1}} = b_{{1}}t + b_{{2}}t^{{2}} , \hspace{2mm} $ we note that the coefficient of $ {t}^{{n + 1}}\hspace{2mm} $on the left side of Eq.  \eqref{EQUATION.4e36cc79-02db-45e8-86b0-cff2382aeaa5} is$ \hspace{2mm}K{ \left( n + 1,{B}_{{2}}^{{n}} \right) } $ and on the right is $ {a}_{{1}}K{ \left( 1,{B}_{{2}}^{{1}} \right) }. $ Therefore we want

       \begin{equation}\label{EQUATION.80d0a2bb-fd20-42e1-a492-671c8154300a}
       K \left( n + 1,{B}_{{2}}^{{n}} \right)  = {a}_{{1}}K{ \left( 1,{B}_{{2}}^{{1}} \right) }.
       \end{equation}

        \noindent Similarly, equating the coefficients of $ {t}^{{n + 2}} $ with a solution $ {B}_{{3}}^{{1}}\hspace{2mm} $yields

       \begin{equation}\label{EQUATION.ac663f3d-5147-46cf-acca-50b3fce49e18}
       K \left( n + 2,{B}_{{3}}^{{n}} \right)  = {a}_{{2}}K{ \left( 2,{B}_{{3}}^{{2}} \right)  + {a}_{{1}}K{ \left( 2,{B}_{{3}}^{{1}} \right) }.}
       \end{equation}

        \noindent When we equate the coefficients of $ t^{{n + 3}}\hspace{2mm} $a similar expression arises with a third term involving the coefficient $ {a}_{{3}} $. In general, to equate the coefficients of $ {t}^{{n + q}} $ on both sides of Eq.  \eqref{EQUATION.4e36cc79-02db-45e8-86b0-cff2382aeaa5} (with a candidate solution $ {B}_{{q + 1}}^{{1}} $) one needs:

       \begin{equation}\label{EQUATION.26e6c154-9a63-44c9-ae45-15c4a1a32ad6}
       K{ \left( n + q,{B}_{{q + 1}}^{{n}} \right) } = \underset{{m = 1}}{\overset{{min{ \left( n - 1,q \right) }}}{ \sum }}{a}_{{m}}K{ \left( q,{B}_{{q + 1}}^{{m}} \right) }\hspace{2mm}q = 1,2,...
       \end{equation}

        \noindent Equation  \eqref{EQUATION.0269eb99-aa8d-4ad2-921d-7baf70e590b2} shows that for any $ s \ge 1 $

       \begin{equation}\label{}
       K{ \left( n + q , {B}_{{q + 1}}^{{n}} \right) } = K{ \left( n + q , {B}_{{q + s}}^{{n}} \right) }
       \end{equation}

       \begin{equation}\label{}
       K{ \left( q , {B}_{{q + 1}}^{{m}} \right) } = K{ \left( q , {B}_{{q - m + 1}}^{{m}} \right) } = K{ \left( q , {B}_{{q + s}}^{{m}} \right) }.
       \end{equation}

        \noindent These equations show that  \eqref{EQUATION.26e6c154-9a63-44c9-ae45-15c4a1a32ad6} implies for any $ s $

       \begin{equation}\label{}
       K{ \left( n + q,{B}_{{q + s}}^{{n}} \right) } = \underset{{m = 1}}{\overset{{min{ \left( n - 1,q \right) }}}{ \sum }}{a}_{{m}}K{ \left( q,{B}_{{q + s}}^{{m}} \right) }\hspace{2mm}.
       \end{equation}

        \noindent Therefore if  \eqref{EQUATION.26e6c154-9a63-44c9-ae45-15c4a1a32ad6} is satisfied, the coefficients of $ t^{{n + q}}\hspace{2mm} $on both sides of Eq.  \eqref{EQUATION.4e36cc79-02db-45e8-86b0-cff2382aeaa5} are also equal when any number of terms $ b_{{k}}t^{{k}}\hspace{2mm} $are added to the partial sum $ {B}_{{q + 1}}^{{1}}. $

       We next write the $ K{ \left( a,{B}_{{u}}^{{v}} \right)  ' s}\hspace{2mm} $appearing in Eq.  \eqref{EQUATION.26e6c154-9a63-44c9-ae45-15c4a1a32ad6} in such a way that the indices $ a $ and $ u $ are equal. Equation  \eqref{EQUATION.0269eb99-aa8d-4ad2-921d-7baf70e590b2} shows that Eq.  \eqref{EQUATION.26e6c154-9a63-44c9-ae45-15c4a1a32ad6} is equivalent to

       \EmptyParagraph

       \begin{equation}\label{EQUATION.84295b17-1eab-4f70-8cbc-cf3003b94ee3}
       K{{ \left( n + q,{B}_{{n + q}}^{{n}} \right) } = \underset{{m = 1}}{\overset{{min{ \left( n - 1,q \right) }}}{ \sum }}{a}_{{m}}K{ \left( q,{B}_{{q}}^{{m}} \right) }\hspace{2mm}q = 1,2...}
       \end{equation}

       \EmptyParagraph

       Given that $ b_{{q}} = K{ \left( q,{B}_{{q}}^{{1}} \right) } , \hspace{2mm} $the task will be to find$ \hspace{2mm}n\hspace{2mm} $ sequences of $ b_{{m}} ' s\hspace{2mm} $for which the infinite system  \eqref{EQUATION.84295b17-1eab-4f70-8cbc-cf3003b94ee3} is satisfied for all $ q.\hspace{2mm} $We will achieve this through a simple discrete dynamical system that generates the $ b_{{q}} ' s $. A change of variable on the dynamical system will show that the $ n\hspace{2mm} $sequences $ { \left\{ b_{{m}} \right\} } $ are of the form $ { \left\{ { \beta }_{{m}}e^{{2k \pi m \times i / n}} \right\} } $ required in Eq.  \eqref{EQUATION.e4267f21-c149-484a-995b-9d3574d58deb} (with $  \left. k = 0 , 1 , ... , n - 1 \right) . $ Finally we will give the convergence conditions for the series $  \sum { \left| { \beta }_{{m}}t^{{m}} \right| } $ and show that the $ x{ \left( t , k \right) } ' s\hspace{2mm} $are the $ n\hspace{2mm} $roots$ \hspace{2mm} $in the sense of  \eqref{EQUATION.d182e62d-c715-4015-8a4a-f408167c8ca0}.

      \section{Matrix formulation}\label{SECTION.425038e3-ade0-4dd3-9c60-be784a885762}

        \noindent We first define the sequence of $ n - 1\hspace{2mm} $dimensional vectors

       \begin{equation}\label{EQUATION.0dd8d22d-af88-48ce-8362-d10640cada3c}
       V{{ \left( q \right) }\overset{{def.}}{{ = }}{ \left( \begin{matrix}K{ \left( q,{B}_{{q}}^{{1}} \right) } & K{ \left( q + 1,{B}_{{q + 1}}^{{2}} \right) } & ... & K{ \left( q + n - 2,{B}_{{q + n - 2}}^{{n - 1}} \right) }\end{matrix} \right)  ' }\hspace{2mm}\hspace{2mm}q = 1,2...}
       \end{equation}

        \noindent where the apostrophe means "transpose" so that $ V{ \left( q \right) } $ is a column vector. We note that

       \begin{equation}\label{EQUATION.cafc54e4-d3ed-4ad4-9956-90cd91d4f0b5}
       V{{ \left( 1 \right) } = { \left( \begin{matrix}{b}_{{1}} & {b}_{{1}}^{{2}} & {b}_{{1}}^{{3}} & ... & {b}_{{1}}^{{n - 1}}\end{matrix} \right)  ' }.}
       \end{equation}

        \noindent Also the $ s - th\hspace{2mm} $component $ V{{ \left( 2 \right) }}_{{s}} $ of $ V{ \left( 2 \right) } $ is$ \hspace{2mm}K{ \left( s + 1 , {B}_{{s + 1}}^{{s}} \right) } $ $  = {{s.b}}_{{2}}{b}_{{1}}^{{s - 1}} $ with $ b_{{2}}\hspace{2mm} $to be determined from Eq.  \eqref{EQUATION.84295b17-1eab-4f70-8cbc-cf3003b94ee3} for $ q = 1 :  $

       \begin{equation}\label{}
       n.b_{{2}}{b}_{{1}}^{{n - 1}} = K{ \left( n + 1 , {B}_{{n + 1}}^{{n}} \right) } = {a}_{{1}}b_{{1}}
       \end{equation}

        \noindent from which $ b_{{2}} = a_{{1}} / { \left( n{b}_{{1}}^{{n - 2}} \right) } $ and therefore

       \begin{equation}\label{EQUATION.c4b29494-056a-49fc-b823-30e68a89fd02}
       V{ \left( 2 \right)  = { \left( \begin{matrix}{a}_{{1}} / n{b}_{{1}}^{{n - 2}} & 2{a}_{{1}} / n{b}_{{1}}^{{n - 3}} & ... & ... & { \left( n - 1 \right) {a}_{{1}}} / n{b}_{{1}}^{{0}}\end{matrix} \right)  ' }.}
       \end{equation}

       \EmptyParagraph

       We next define the sequence of $ n - 1\hspace{2mm} $dimensional vectors $ X{ \left( q \right) } $ obtained as the convolution of the $ {b}_{{i}} ' s\hspace{2mm} $with the vectors $ V{ \left( q \right)  : } $

       \begin{equation}\label{EQUATION.b57e840f-9395-475e-a584-b00621c935cf}
       X{{ \left( q - 1 \right) }\overset{{def.}}{{ = }}\underset{{m = 2}}{\overset{{q - 1}}{ \sum }}b_{{m}}V{ \left( q - m + 1 \right) }\hspace{2mm}q = 3,4,....\hspace{2mm}}
       \end{equation}

        \noindent We let $ V{{ \left( q \right) }}_{{r}} $ and $ X{{ \left( q \right) }}_{{r}} $ denote the $ r - th\hspace{2mm} $ components of $ V{ \left( q \right) }\hspace{2mm} $and $ X{ \left( q \right) }. $ The left-hand side of Eq.  \eqref{EQUATION.453f5f3b-86ea-4183-9ce2-6adf54a01a37} is $ V{{{ \left( q \right) }}}_{{s + 1}} $ and Eq.  \eqref{EQUATION.453f5f3b-86ea-4183-9ce2-6adf54a01a37} with $ s = 1\hspace{2mm} $is then

       \begin{equation}\label{EQUATION.4419766a-2c29-4c61-98ee-bcd0e8c272e7}
       V{{ \left( q \right) }}_{{2}} = b_{{1}}V{{ \left( q \right) }}_{{1}} + X{{ \left( q - 1 \right) }}_{{1}} + {{V{ \left( q \right) }}}_{{1}}{b}_{{1}}^{{}} = 2b_{{1}}V{{ \left( q \right) }}_{{1}} + X{{ \left( q - 1 \right) }}_{{1}}.
       \end{equation}

        \noindent With $ s = 2\hspace{2mm} $ Eq.  \eqref{EQUATION.453f5f3b-86ea-4183-9ce2-6adf54a01a37} is

       \begin{equation}\label{}
       V{{ \left( q \right) }}_{{3}} = b_{{1}}V{{ \left( q \right) }}_{{2}} + X{{ \left( q - 1 \right) }}_{{2}} + {{V{ \left( q \right) }}}_{{1}}{b}_{{1}}^{{2}}
       \end{equation}

        \noindent and more generally

       \begin{equation}\label{EQUATION.976a3359-086f-432c-b205-83610ab76d27}
       V{{ \left( q \right) }}_{{s + 1}} = b_{{1}}V{{ \left( q \right) }}_{{s}} + X{{ \left( q - 1 \right) }}_{{s}} + {{V{ \left( q \right) }}}_{{1}}{b}_{{1}}^{{s}},\hspace{2mm}s = 1,2,...,n - 2.
       \end{equation}

        \noindent Equation  \eqref{EQUATION.453f5f3b-86ea-4183-9ce2-6adf54a01a37} with $ s = n - 1 $ is

       \begin{equation}\label{EQUATION.b655b07f-5609-4115-8fd0-97ad9827f566}
        \left. K{ \left( q + n - 1 \right. },{B}_{{q + n - 1}}^{{n}} \right)  = b_{{1}}V{{ \left( q \right) }}_{{n - 1}} + X{{ \left( q - 1 \right) }}_{{n - 1}} + {{V{ \left( q \right) }}}_{{1}}{b}_{{1}}^{{n - 1}}.
       \end{equation}

        \noindent If we define the $ n - 1\hspace{2mm} $ dimensional square matrix

       \begin{equation}\label{}
       M_{{1}} = { \left( \begin{matrix}2{b}_{{1}} &  - 1 & 0 & 0 & ... & 0\\
            {b}_{{1}}^{{2}} & {b}_{{1}} &  - 1 & 0 & ... & 0\\
            {b}_{{1}}^{{3}} & 0 & {b}_{{1}} &  - 1 & ... & 0\\
            ... & ... & ... & ... & ... & ...\\
            {b}_{{1}}^{{n - 2}} & 0 & ... & 0 & {b}_{{1}} &  - 1\\
            {b}_{{1}}^{{n - 1}} & 0 & ... & ... & ... & {b}_{{1}}\end{matrix} \right) }
       \end{equation}

        \noindent then Eqs.  \eqref{EQUATION.4419766a-2c29-4c61-98ee-bcd0e8c272e7}- \eqref{EQUATION.b655b07f-5609-4115-8fd0-97ad9827f566} (for $  \left. s = 1,2,...,n - 1 \right) \hspace{2mm} $can be written compactly as

       \begin{equation}\label{EQUATION.41c1f466-f964-4583-8bb9-c3a2e50921c2}
        - X{{ \left( q - 1 \right) } = M_{{1}}V{{ \left( q \right) } - { \left( \begin{matrix}0\\
            0\\
            ...\\
            ...\\
            K{ \left( q + n - 1 , {B}_{{q + n - 1}}^{{n}} \right) }\end{matrix} \right) }}}
       \end{equation}

        \noindent The goal is to express $ V{ \left( q \right) } $ as a function of past vectors $ V{ \left( q - 1 \right) },\hspace{2mm}V{ \left( q - 2 \right) },...,V{ \left( 1 \right) }.\hspace{2mm} $The inverse of $ {M}_{{1}} $ is

       \begin{equation}\label{EQUATION.13b2bd16-6f87-4f3d-a744-ecef9a7acb63}
       {M}_{{1}}^{{ - 1}} = \frac {{1}}{{n{b}_{{1}}^{{n - 1}}}}{ \left( \begin{matrix}{b}_{{1}}^{{n - 2}} & {b}_{{1}}^{{n - 1}} & ... & {b}_{{1}}^{{2}} & {b}_{{1}} & 1\\
            {{ - { \left( n - 2 \right) }b}}_{{1}}^{{n - 1}} & {{2b}}_{{1}}^{{n - 2}} & {{2b}}_{{1}}^{{n - 3}} & ... & {{2b}}_{{1}}^{{2}} & {{2b}}_{{1}}^{{}}\\
             - { \left( n - 3 \right) }{b}_{{1}}^{{n}} & {{ - { \left( n - 3 \right) }b}}_{{1}}^{{n - 1}} & 3{b}_{{1}}^{{n - 2}} & ... & 3{b}_{{1}}^{{3}} & 3{b}_{{1}}^{{2}}\\
            ... & ... & ... & ... & ... & ...\\
            {{ - 2b}}_{{1}}^{{2n - 5}} &  - 2{b}_{{1}}^{{2n - 6}} & ... &  - 2{b}_{{1}}^{{n - 1}} & { \left( n - 2 \right) {b}_{{1}}^{{n - 2}}} & { \left( n - 2 \right) }{b}_{{1}}^{{n - 3}}\\
            {{ - b}}_{{1}}^{{2n - 4}} & {{ - b}}_{{1}}^{{2n - 5}} & {{ - b}}_{{1}}^{{2n - 6}} & ... &  - {b}_{{1}}^{{n - 1}} & { \left( n - 1 \right) }{b}_{{1}}^{{n - 2}}\end{matrix} \right) }
       \end{equation}

        \noindent If $ {{ \left( {M}_{{1}}^{{ - 1}} \right) }}_{{i,j}} $ is the entry in the i-th row, j-th entry of this inverse, Eq.  \eqref{EQUATION.13b2bd16-6f87-4f3d-a744-ecef9a7acb63} can be written as

       \begin{equation}\label{}
       {{ \left( {M}_{{1}}^{{ - 1}} \right) }}_{{i,j}} = \frac {{1}}{{n{b}_{{1}}^{{n - 1}}}}{ \left\{ \begin{matrix} - { \left( n - i \right) {b}_{{1}}^{{n + i - j - 2}}\hspace{2mm}if\hspace{2mm}i > j}\\
            i{b}_{{1}}^{{n + i - j - 2}}\hspace{2mm}if\hspace{2mm}i \le j\end{matrix} \right. }
       \end{equation}

       \EmptyParagraph

       We multiply both sides of Eq.  \eqref{EQUATION.41c1f466-f964-4583-8bb9-c3a2e50921c2} by $ {M}_{{1}}^{{ - 1}} $ to obtain

       \begin{equation}\label{EQUATION.e3efea2f-dcc6-4b88-91cd-68f3d61d6207}
       V{{ \left( q \right) } =  - {M}_{{1}}^{{ - 1}}X{{ \left( q - 1 \right) } + {M}_{{1}}^{{ - 1}} \left( \begin{matrix}0\\
            0\\
            ...\\
            ...\\
            {K{ \left( q + n - 1 , {B}_{{q + n - 1}}^{{n}} \right) }}\end{matrix} \right) }.}
       \end{equation}

        \noindent From Eq.  \eqref{EQUATION.84295b17-1eab-4f70-8cbc-cf3003b94ee3} we note that

       \begin{equation}\label{}
       K{ \left( q + n - 1 , {B}_{{q + n - 1}}^{{n}} \right)  = \underset{{m = 1}}{\overset{{min{ \left( n - 1,q - 1 \right) }}}{ \sum }}{a}_{{m}}K{ \left( q - 1,{B}_{{q - 1}}^{{m}} \right) }}
       \end{equation}

       \begin{equation}\label{}
        = \underset{{m = 1}}{\overset{{min{ \left( n - 1,q - 1 \right) }}}{ \sum }}{a}_{{m}}V{{ \left( q - m \right) }}_{{m}} = \underset{{m = 1}}{\overset{{min{ \left( n - 1,q - 1 \right) }}}{{ \sum }}}{{{{ \left[ A_{{m}}.V{ \left( q - m \right) } \right] }}}}_{{}}
       \end{equation}

        \noindent where each $ {A}_{{p}}\hspace{2mm}{ \left( p = 1,2,... , n - 1 \right) } $ is an $ { \left( n - 1 \right) } $-dimensional row vector with $ a_{{p}}\hspace{2mm} $in $ p - th $ position and zeros elsewhere. (Here and elsewhere an expression in square brackets will generally represent a scalar product).

       The vector $ V{ \left( q \right) } $ of Eq.  \eqref{EQUATION.e3efea2f-dcc6-4b88-91cd-68f3d61d6207} is now

       \begin{equation}\label{EQUATION.91aaf7ef-2004-4f50-aca2-74c9f1c24b3c}
       V{ \left( q \right) } =  - {M}_{{1}}^{{ - 1}}X{{ \left( q - 1 \right) } + \frac {{\underset{{m = 1}}{\overset{{min{ \left( n - 1,q - 1 \right) }}}{{ \sum }}} \left[ {A}_{{m}}.V{ \left( q - m \right) } \right] }}{{n{b}_{{1}}^{{n - 1}}}}{ \left( \begin{matrix}1\\
            2{b}_{{1}}\\
            3{b}_{{1}}^{{2}}\\
            ...\\
            { \left( {n - 1} \right) {b}_{{1}}^{{n - 2}}}\end{matrix} \right) }.}
       \end{equation}

       \EmptyParagraph

       \EmptyParagraph

       For each one of the n values $ \hspace{2mm}{{ \rho }}^{{1 / n}}{e}^{{i{ \left(  \theta  + 2k \pi  \right)  / n}}}\hspace{2mm} $ ($  \left. k = 0,1,...,n - 1 \right) \hspace{2mm} $ of $ b_{{1}}\hspace{2mm} $the sequence $ V{ \left( q \right) } $ yields all the subsequent $ {b}_{{i}} ' s\hspace{2mm} $since $ {b}_{{i}} $ is the first component of $ V{ \left( i \right) }.\hspace{2mm} $We thus have the $ n\hspace{2mm} $sequences $ { \left\{ {b}_{{k}} \right\} }\hspace{2mm} $that will provide the $ n $ roots (under the right as-yet-unproved convergence conditions).

       A change of variable on the sequence $ V{ \left( m \right) }\hspace{2mm} $which defines a new sequence $ W{ \left( m \right) } $ will now show that the $ n $ values of each $ b_{{m\hspace{2mm}}} $can be expressed in the form $ { \beta }_{{m}}e^{{2k \pi m \times i / n}} $ $ { \left( k = 0 , 1 , ... , n - 1 \right) } $ where a single sequence $ { \left\{ { \beta }_{{m}} \right\} }\hspace{2mm} $is defined recursively as a function of the $ a_{{k}} ' s\hspace{2mm} $.

       \EmptyParagraph

       We will first need the diagonal matrix

       \begin{equation}\label{}
       B\overset{{def.}}{{ = }} \left( \begin{matrix}{b}_{{1}}^{{n - 1}} & 0 & 0 & 0 & 0\\
            0 & {b}_{{1}}^{{n - 2}} & 0 & 0 & 0\\
            0 & 0 & {b}_{{1}}^{{n - 3}} & 0 & 0\\
            ... & ... & ... & ... & 0\\
            ... & ... & ... & 0 & {b}_{{1}}^{{1}}\end{matrix} \right) 
       \end{equation}

        \noindent and

       \begin{equation}\label{EQUATION.b70e0749-3e15-4609-9d3f-fa3a5aad21ca}
       M\overset{{def.}}{{ = }}{\frac {{1}}{{n}} \left( \begin{matrix}1 & 1 & ... & ... & 1\\
            { -  \left( n - 2 \right) } & 2 & 2 & ... & 2\\
             - { \left( n - 3 \right) } &  - { \left( n - 3 \right) } & 3 & ... & 3\\
            ... & ... & ... & ... & ...\\
             - 1 &  - 1 & ... &  - 1 & n - 1\end{matrix} \right) }
       \end{equation}

       \EmptyParagraph

       The matrix $ M $ is the inverse $ {M}_{{1}}^{{ - 1}} $ of  \eqref{EQUATION.13b2bd16-6f87-4f3d-a744-ecef9a7acb63} in which $ {b}_{{1}} $ would be set to 1. With these notations

       \begin{equation}\label{}
       {M}_{{1}}^{{ - 1}} = \frac {{1}}{{{b}_{{1}}}}{B}^{{ - 1}}M\hspace{2mm}B
       \end{equation}

        \noindent The sequence $ W{ \left( q \right) }\hspace{2mm} $derived from $ V{ \left( q \right) }\hspace{2mm} $is defined as

       \begin{equation}\label{EQUATION.ca293d2d-3eb8-4031-a264-b8e2983bf9be}
       W{{ \left( q \right) }\overset{{def.}}{{ = }}{b}_{{1}}^{{1 - q}}B.V{ \left( q \right) },\hspace{2mm}q = 1,2,....}
       \end{equation}

       \EmptyParagraph

       \EmptyParagraph

       We will see that the sequence$ \hspace{2mm}W{ \left( q \right) } $ is independent of the particular $ \hspace{2mm}b_{{1}}\hspace{2mm} $chosen among its $ n\hspace{2mm} $possible values $ { \rho }^{{1 / n}}e^{{i{ \left(  \theta  + 2k \pi  \right)  / n}}} $ ($  \left. k = 0 , 1 , ... , n - 1 \right) . $ The fact that each one of the $ b_{{m}} ' s $ is of the form $ { \beta }_{{m}}e^{{2k \pi m \times i / n}} $$ { \left( k = 0 , 1 , ... , n - 1 \right) } $ follows then from the fact that the first components $ V{{ \left( q \right) }}_{{1\hspace{2mm}}}and\hspace{2mm}W{{ \left( q \right) }}_{{1}} $ satisfy

       \begin{equation}\label{EQUATION.a72906b5-d4e2-4d7a-abef-7c6a88f55f88}
       V{{ \left( q \right) }}_{{1}} = W{{ \left( q \right) }}_{{1}}{b}_{{1}}^{{q}} / a_{{0}}.
       \end{equation}

       \EmptyParagraph

       In the theorem below we give for the sequence $ W{ \left( q \right) } $ an expression analogous to  \eqref{EQUATION.91aaf7ef-2004-4f50-aca2-74c9f1c24b3c} for $ V{ \left( q \right) }\hspace{2mm} $and provide the expressions for the $ { \beta }_{{m}} ' s. $$  $
     
       \begin{theorem}\label{THEOREM.b219007d-b5c1-4f0e-9c17-ec591b878ec2}
     
        \noindent The sequence W(q) is defined as

       \begin{equation}\label{EQUATION.663ecc6d-ec42-4137-ba7b-478afd60c800}
       W{ \left( 1 \right) } = a_{{0}}{ \left( \begin{matrix}1 & 1 & ... & ... & 1\end{matrix} \right) } ' 
       \end{equation}

       \begin{equation}\label{EQUATION.1f610449-1951-4969-a6b8-62605182185a}
       W{ \left( 2 \right) } = \frac {{a_{{1}}}}{{n}}{ \left( \begin{matrix}1 & 2 & ... & ... & n - 1\end{matrix} \right) } ' 
       \end{equation}

       \begin{equation}\label{EQUATION.4ee97e62-5b23-4a08-b90e-19f6a219dc0e}
       W{{ \left( q \right) } = \frac {{1}}{{a_{{0}}}}{ \left[  - M\underset{{p = 1}}{\overset{{q - 2}}{ \sum }}{ \left[ u.W{ \left( p + 1 \right) } \right] }W{ \left( q - p \right) } + U\underset{{p = 1}}{\overset{{min{ \left( n - 1 , q - 1 \right) }}}{ \sum }}{ \left[ A_{{p}}.W{ \left( q - p \right) } \right] } \right]  , q = 3 , 4 , ...}}
       \end{equation}

        \noindent where:

       1. u is the (n-1)-dimensional row vector having 1 in first position and zeros elsewhere.

       2. U is the (n-1)-dimensional column vector $ \frac {{1}}{{n}}{ \left( \begin{matrix}1 & 2 & ... & ... & n - 1\end{matrix} \right)  ' }. $

       The $ {{ \beta }}_{{m}} ' s $ that define the $ n $ sequences

       \begin{equation}\label{EQUATION.fb2e134f-95a6-4fa2-b570-b69c33a2e615}
       b_{{m}} = { \beta }_{{m}}e^{{2k \pi m \times i / n}} , k = 0 , 1 , ... , n - 1 ; m = 1 , 2 , ...
       \end{equation}

        \noindent are

       \begin{equation}\label{EQUATION.3e6e9933-f6ec-48a5-9d58-d9402b658d2e}
       { \beta }_{{m}} = { \rho }^{{m / n}} \times e^{{i{ \left( m \theta  \right)  / n}}} \times W{{ \left( m \right) }}_{{1}} / a_{{0}} , \hspace{2mm}m = 1 , 2 , ...
       \end{equation}

       \EmptyParagraph
     
    \end{theorem}\begin{proof}
     
        \noindent The expressions of  \eqref{EQUATION.663ecc6d-ec42-4137-ba7b-478afd60c800}- \eqref{EQUATION.1f610449-1951-4969-a6b8-62605182185a} follow from the definition in Eq.  \eqref{EQUATION.ca293d2d-3eb8-4031-a264-b8e2983bf9be} and Eqs.  \eqref{EQUATION.cafc54e4-d3ed-4ad4-9956-90cd91d4f0b5} and  \eqref{EQUATION.c4b29494-056a-49fc-b823-30e68a89fd02} which define $ V{ \left( 1 \right) }\hspace{2mm} $and $ V{ \left( 2 \right) }. $

       Using Eq.  \eqref{EQUATION.b57e840f-9395-475e-a584-b00621c935cf} to express $ X{ \left( q - 1 \right) }\hspace{2mm} $ appearing in $ V{ \left( q \right) } $ of Eq.  \eqref{EQUATION.91aaf7ef-2004-4f50-aca2-74c9f1c24b3c} yields

       \EmptyParagraph

       \[
       V{ \left( q \right) } = {b}_{{1}}^{{q - 1}}B^{{ - 1}}W{{ \left( q \right) } =  - {M}_{{1}}^{{ - 1}}{ \left( \underset{{m = 2}}{\overset{{q - 1}}{ \sum }}b_{{m}}{b}_{{1}}^{{q - m}}B^{{ - 1}}W{ \left( q - m + 1 \right) } \right) } + }
       \]

       \begin{equation}\label{EQUATION.fe8e0790-fa4c-4873-93d7-b93de7fbbf35}
       \frac {{\underset{{p = 1}}{\overset{{min{ \left( n - 1 , q - 1 \right) }}}{ \sum }}{ \left[ A_{{p}}.{ \left( {{{b}_{{1}}^{{q - p - 1}}B}}^{{ - 1}}W{ \left( q - p \right) } \right) } \right] }}}{{n{b}_{{1}}^{{n - 1}}}}{ \left( \begin{matrix}1\\
            2b_{{1}}\\
            3{b}_{{1}}^{{2}}\\
            ...\\
            { \left( n - 1 \right) {b}_{{1}}^{{n - 2}}}\end{matrix} \right) } , \hspace{2mm}q = 3 , 4 , ...
       \end{equation}

        \noindent Bearing in mind that:

       $  \ast \hspace{2mm}{{B \times M}}_{{1}}^{{ - 1}}{{ \times B}}^{{ - 1}} = M / b_{{1}} $; $ {\hspace{2mm}b}_{{1}}^{{n}} = a_{{0}} $

       $  \ast \hspace{2mm} $each $ b_{{m}} $ is $ {b}_{{1}}^{{m - n}} $ multiplied by the first component $ { \left[ u.W{ \left( m \right) } \right] } $ of $ W{ \left( m \right) } $

       $  \ast \hspace{2mm} $the $ p - th $ component of $ B^{{ - 1}}W{ \left( q - p \right) } $ is $ {b}_{{1}}^{{p - n}}W{{ \left( q - p \right) }}_{{p}} = {b}_{{1}}^{{p}}W{{ \left( q - p \right) }}_{{p}} / a_{{0}} ,  $

       we multiply both sides of Eq.  \eqref{EQUATION.fe8e0790-fa4c-4873-93d7-b93de7fbbf35} by $ {b}_{{1}}^{{1 - q}}B $ to obtain

       \begin{equation}\label{}
       W{ \left( q \right)  = }\frac {{1}}{{a_{{0}}}}{ \left[  - M\underset{{m = 2}}{\overset{{q - 1}}{ \sum }}{{ \left[ u.W{ \left( m \right) } \right] }W{ \left( q - m + 1 \right) } + U\underset{{p = 1}}{\overset{{min{ \left( n - 1 , q - 1 \right) }}}{ \sum }}{ \left[ A_{{p}}.W{ \left( q - p \right) } \right] }} \right]  , q = 3 , 4 , ...}
       \end{equation}

        \noindent which is Eq.  \eqref{EQUATION.4ee97e62-5b23-4a08-b90e-19f6a219dc0e} once the index $ m $ in the first sum is set equal to $ p + 1 $ with $ p $ going from 1 to $ q - 2. $

       In view of Eq.  \eqref{EQUATION.a72906b5-d4e2-4d7a-abef-7c6a88f55f88} the first component $ b_{{m}}\hspace{2mm} $of $ V{ \left( m \right) } $ is equal to

       \begin{equation}\label{}
       b_{{m}} = {b}_{{1}}^{{m}}W{{ \left( m \right) }}_{{1}} / a_{{0}} = { \rho }^{{m / n}} \times e^{{m{ \left(  \theta  + 2k \pi  \right) }i / n}}W{{ \left( m \right) }}_{{1}} / a_{{0}}
       \end{equation}

        \noindent which is  \eqref{EQUATION.fb2e134f-95a6-4fa2-b570-b69c33a2e615} with $ { \beta }_{{m}}\hspace{2mm} $ given in Eq.  \eqref{EQUATION.3e6e9933-f6ec-48a5-9d58-d9402b658d2e}.
     
       \end{proof}

      \section{Convergence results}\label{SECTION.b7a6c94d-ac14-4555-a497-5aa6412e23e1}

        \noindent In order to establish convergence conditions for the series of  \eqref{EQUATION.c5df2467-e55b-4766-9170-787c0bdffe48} we need to assess the growth of the first components $ W{{ \left( m \right) }}_{{1}}\hspace{2mm} $of the vectors $ W{ \left( m \right) }. $ We will use the $ {\mathrm{ \ell }}_{{ \infty }}\hspace{2mm} $norm $ {{{ \left| V \right| }}}_{{ \infty }} = \underset{{i}}{{max}}{ \left| V_{{i}} \right| } $ of a complex vector $ V{ = { \left( v_{{i}} \right) }.} $ We also define the row-sum norm

       \begin{equation}\label{}
       { \left\| A \right\| } = {max}_{{i}}\hspace{2mm}{ \sum }_{{j = 1}}^{{m}}{ \left| a_{{ij}} | \hspace{2mm} \right. }
       \end{equation}

        \noindent of an $ m -  $dimensional square matrix $ A = { \left( a_{{ij}} \right) }\hspace{2mm} $and recall that $ {{{{ \left| AV \right| }}}_{{ \infty }} \le {{ \left\| A \right\| }{{{{ \left| V \right| }}}_{{ \infty }}.}}} $

       If we define the function of two variables

       \begin{equation}\label{}
       z{ \left( u , n \right) } =  - 2u^{{2}} + u{ \left( 2n + 1 \right) } - n
       \end{equation}

        \noindent then the norm of the $ { \left( n - 1 \right) } -  $dimensional matrix $ M $ of Eq.  \eqref{EQUATION.b70e0749-3e15-4609-9d3f-fa3a5aad21ca} is

       \begin{equation}\label{}
       {{ \left\| M \right\| } = \frac {{1}}{{n}}max{ \left( z{{ \left[ floor{ \left( n / 2 + 1 / 4 \right)  , n} \right] } , z{ \left[ ceil{ \left( n / 2 + 1 / 4 \right)  , n} \right] }}\hspace{2mm} \right) }.}
       \end{equation}

        \noindent where the $ floor{ \left( x \right) } $ and $ ceil{ \left( x \right) }\hspace{2mm} $functions are the largest integer smaller than $ x $ and the smallest integer less than $ x. $

       We next define the sequence of modulii

       \begin{equation}\label{}
       w{ \left( q \right) }\overset{{def.}}{{ = }}{{{{ \left| W{ \left( q \right) } \right| }}}_{{ \infty }}\hspace{2mm} , q = 1 , 2 , ...}
       \end{equation}

        \noindent which are upper bounds for the modulii of the $ W{{ \left( q \right) }}_{{1}} ' s\hspace{2mm} $of interest. We also define the maximum modulus of the coefficients $ a_{{k}}\hspace{2mm} $other than $ a_{{0}} :  $

       \begin{equation}\label{}
        \alpha \overset{{def.}}{{ = }}\underset{{k = 1 , 2 , ... , n - 1}}{{max}}{ \left| a_{{k}} | . \right. }
       \end{equation}

        \noindent We note that $ w{ \left( 1 \right)  = { \left| a_{{0}} | \hspace{2mm} \right. }} $and $ w{ \left( 2 \right)  \le  \alpha .\hspace{2mm}} $From Eq.  \eqref{EQUATION.4ee97e62-5b23-4a08-b90e-19f6a219dc0e} we now have

       \begin{equation}\label{EQUATION.a365d4a5-8923-45e4-a665-1d6c773162f6}
       w{ \left( q \right) } \le \frac {{1}}{{{ \left| a_{{0}} \right| }}}{ \left( \underset{{p = 1}}{\overset{{q - 2}}{ \sum }}w{ \left( p + 1 \right) }{{{ \left\| M \right\| } \times w{ \left( q - p \right) }} +  \alpha \underset{{p = 1}}{\overset{{min{ \left( n - 1 , q - 1 \right) }}}{{ \sum }}}w{ \left( q - p \right) }} \right)  , \hspace{2mm}q = 3 , 4 , ...}
       \end{equation}

        \noindent with an important distinction to be made between the cases $ q \le n $ and $ q > n.\hspace{2mm}   $In the former case the $ p\hspace{2mm} $in the second sum on the right-hand side of  \eqref{EQUATION.a365d4a5-8923-45e4-a665-1d6c773162f6} goes to $ q - 1\hspace{2mm} $with a last term $ w{ \left( 1 \right) } = { \left| a_{{0}} | .\hspace{2mm} \right. } $The inequality of  \eqref{EQUATION.a365d4a5-8923-45e4-a665-1d6c773162f6} is then

       \begin{equation}\label{EQUATION.1c2b7541-d50a-43c5-94dd-63ad5157fa42}
       w{ \left( q \right) } \le \frac {{1}}{{{ \left| a_{{0}} \right| }}}{{ \left( \underset{{p = 1}}{\overset{{q - 2}}{ \sum }}{ \left[ { \left\| M \| w{ \left( p + 1 \right) } +  \alpha  \right. } \right] }w{ \left( q - p \right) } \right) } +  \alpha  , \hspace{2mm}q = 3 , 4 , ... , n}
       \end{equation}

        \noindent We next define

       \begin{equation}\label{}
        \mu  = { \left\| M \|  / { \left| a_{{0}} \right| } \right. }
       \end{equation}

        \noindent and the sequence of functions

       \begin{equation}\label{}
       { \sigma }_{{p}}{ \left(  \alpha  \right) } = { \left\{ \begin{matrix} \alpha  / { \left\| M \right\| }{ , \hspace{2mm}if\hspace{2mm}p \le n}\\
            0 , \hspace{2mm}if\hspace{2mm}p > n\end{matrix} \right\}  , \hspace{2mm}p = 2 , 3 , ...}
       \end{equation}

        \noindent In the case $ q > n $ the index $ p $ in the second sum on the right-hand side of  \eqref{EQUATION.a365d4a5-8923-45e4-a665-1d6c773162f6} goes to $ n - 1\hspace{2mm} $ and  \eqref{EQUATION.a365d4a5-8923-45e4-a665-1d6c773162f6} can now be written

       \begin{equation}\label{}
       w{ \left( q \right) } \le  \mu { \left( \underset{{p = 1}}{\overset{{q - 2}}{ \sum }}{ \left[ w{ \left( p + 1 \right) } + { \sigma }_{{p + 1}}{ \left(  \alpha  \right) } \right] }w{ \left( q - p \right) } \right) }
       \end{equation}

       \begin{equation}\label{EQUATION.196407ef-81f2-4302-aeaf-f602b93f2754}
        \le  \mu { \left( \underset{{p = 1}}{\overset{{q - 2}}{ \sum }}{ \left[ w{ \left( p + 1 \right) } + { \sigma }_{{p + 1}}{ \left(  \alpha  \right) } \right] }{ \left[ w{ \left( q - p \right) } + { \sigma }_{{q - p}}{ \left(  \alpha  \right) } \right] } \right)  , \hspace{2mm}q = n + 1 , n + 2 , ...}
       \end{equation}

        \noindent where the purpose of this last inequality is to bound the sequence $ { \left\{ w{ \left( q \right) } \right\} }\hspace{2mm} $by a particular type of discrete convolution which will be considered below.

       We first define recursively a sequence of nonnegative functions $ S_{{q}}{ \left(  \alpha  , { \left| a_{{0}} \right| } \right) } :  $

       \begin{equation}\label{EQUATION.d6aaf9f6-1cf5-4f12-8f65-70cb8593d6dc}
       S_{{2}}{ \left(  \alpha  , { \left| a_{{0}} \right| } \right) }\overset{{def.}}{{ = }} \alpha  + { \sigma }_{{2}}{ \left(  \alpha  \right) }
       \end{equation}

       \begin{equation}\label{EQUATION.8a0b6a92-3b23-4597-b3ff-658733e697f8}
       S_{{q}}{ \left(  \alpha  , { \left| a_{{0}} \right| } \right) }\overset{{def.}}{{ = }}\frac {{1}}{{{{ \left| a_{{0}} \right| }}}}{ \left( \underset{{p = 1}}{\overset{{q - 2}}{ \sum }}{ \left[ { \left\| M \| S_{{p + 1}} \left(  \alpha  , { \left| a_{{0}} \right| } \right)  +  \alpha  \right. } \right] }S_{{q - p}}{ \left(  \alpha  , { \left| a_{{0}} \right| } \right) } \right)  +  \alpha  + { \sigma }_{{q}}{ \left(  \alpha  \right)  , q = 3 , 4 , ... , n}}
       \end{equation}

        \noindent where the functional notation emphasizes for future reference a dependence of $ S_{{q}}\hspace{2mm} $on $  \alpha \hspace{2mm} $and $ { \left| a_{{0}} | \hspace{2mm} \right. } $(even though $ S_{{2}}{ \left(  \alpha  , { \left| a_{{0}} \right| } \right) }\hspace{2mm} $of  \eqref{EQUATION.d6aaf9f6-1cf5-4f12-8f65-70cb8593d6dc} does not depend on $ { \left| a_{{0}} \right| } $).

       We next continue the sequence of $ S_{{q}} ' s $ with

       \begin{equation}\label{EQUATION.0032ecb7-5a69-4af0-80d5-1cce7f62fcb4}
       S_{{q}}{ \left(  \alpha  , { \left| a_{{0}} \right| } \right) }\overset{{def.}}{{ = }} \mu { \left( \underset{{p = 1}}{\overset{{q - 2}}{ \sum }}S_{{p + 1}}{ \left(  \alpha  , { \left| a_{{0}} \right| } \right) }S_{{q - p}}{ \left(  \alpha  , { \left| a_{{0}} \right| } \right) } \right) } , q = n + 1 , n + 2 , ...
       \end{equation}
       \begin{lemma}\label{LEMMA.66642c0c-ec86-48b4-a505-ccd6eab1ef85}
     
        \noindent With the notations given above,

       \begin{equation}\label{EQUATION.e89e2e3b-0afd-40bb-8df9-8190e5f926b8}
       w{ \left( q \right) } \le w{ \left( q \right) } + { \sigma }_{{q}}{ \left(  \alpha  \right) } \le S_{{q}}{ \left(  \alpha  , { \left| a_{{0}} \right| } \right) } , \hspace{2mm}q = 2 , 3 , ...
       \end{equation}

       \EmptyParagraph
     
    \end{lemma}\begin{proof}
     
        \noindent The inequality of  \eqref{EQUATION.e89e2e3b-0afd-40bb-8df9-8190e5f926b8} is true for $ q = 2\hspace{2mm} $because $ w{ \left( 2 \right) } \le  \alpha . $ In view of  \eqref{EQUATION.1c2b7541-d50a-43c5-94dd-63ad5157fa42} and  \eqref{EQUATION.8a0b6a92-3b23-4597-b3ff-658733e697f8} it is then also true for $ q = 3 , 4 , ... , n $. For $ q = n + 1 , \hspace{2mm} $  \eqref{EQUATION.196407ef-81f2-4302-aeaf-f602b93f2754} yields

       \begin{equation}\label{EQUATION.76fd0257-2e11-48ba-813f-870afec65ee0}
       w{ \left( n + 1 \right) } = w{ \left( n + 1 \right) } + { \sigma }_{{n + 1}}{{ \left(  \alpha  \right) } \le  \mu { \left( \underset{{p = 1}}{\overset{{n - 1}}{ \sum }}{ \left[ w{ \left( p + 1 \right) } + { \sigma }_{{p + 1}}{ \left(  \alpha  \right) } \right] { \left[ w{ \left( n + 1 - p \right) } + { \sigma }_{{n + 1 - p}}{ \left(  \alpha  \right) } \right] }} \right) }}
       \end{equation}

        \noindent $  \le  \mu {{ \left( \underset{{p = 1}}{\overset{{n - 1}}{ \sum }}{S_{{p + 1}}{ \left(  \alpha  , { \left| a_{{0}} \right| } \right) } \times S_{{n + 1 - p}}{ \left(  \alpha  , { \left| a_{{0}} \right| } \right) }} \right)  = S_{{n + 1}}{ \left(  \alpha  , { \left| a_{{0}} \right| } \right) }.}} $

       To prove the result by induction we assume  \eqref{EQUATION.e89e2e3b-0afd-40bb-8df9-8190e5f926b8} is true up to order $ n + r. $ Then at order $ n + r + 1 $ we have

       \begin{equation}\label{}
       w{ \left( n + r + 1 \right) } = w{ \left( n + r + 1 \right) } + { \sigma }_{{n + r + 1}}{ \left(  \alpha  \right) }
       \end{equation}

       \begin{equation}\label{}
        \le  \mu { \left( \underset{{p = 1}}{\overset{{n + r - 1}}{ \sum }}{ \left[ w{ \left( p + 1 \right) } + { \sigma }_{{p + 1}}{ \left(  \alpha  \right) } \right] { \left[ w{ \left( n + r + 1 - p \right) } + { \sigma }_{{n + r + 1 - p}}{ \left(  \alpha  \right) } \right] }} \right) }
       \end{equation}

       \begin{equation}\label{}
        \le  \mu  \left( \underset{{p = 1}}{\overset{{n + r - 1}}{ \sum }}{S_{{p + 1}}{ \left(  \alpha  , { \left| a_{{0}} \right| } \right) } \times S_{{n + r + 1 - p}}{ \left(  \alpha  , { \left| a_{{0}} \right| } \right) }} \right)  = S_{{n + r + 1}}{ \left(  \alpha  , { \left| a_{{0}} \right| } \right) }
       \end{equation}

        \noindent which completes the proof.

       \EmptyParagraph
     
    \end{proof}
    
       \EmptyParagraph

       We now provide some definitions and results pertaining to convolution-type sequences such as the $ S_{{q}} ' s $ of  \eqref{EQUATION.0032ecb7-5a69-4af0-80d5-1cce7f62fcb4}. 
     
       \begin{proposition}\label{PROPOSITION.bd64019d-cd55-42a1-81d9-61abb499b4ed}
     
        \noindent For any scalar $  \mu  , \hspace{2mm} $a $  \mu  -  $convolution of order $ m $ is an infinite sequence $ { \left\{ u_{{k}} \right\} }\hspace{2mm} $consisting of $ m $ initial scalars $ u_{{1}} , \hspace{2mm}u_{{2}} , .... , u_{{m}} $ with subsequent terms defined as

       \begin{equation}\label{}
       u_{{q}} =  \mu { \left( \underset{{p = 1}}{\overset{{q - 1}}{ \sum }}u_{{p}}u_{{q - p}} \right) }\hspace{2mm} , \hspace{2mm}q = m + 1 , m + 2 , ...
       \end{equation}

        \noindent For a $  \mu  -  $convolution $ { \left\{ v_{{q}} \right\} } $ of order $ m = 1 , \hspace{2mm} $we have

       \begin{equation}\label{EQUATION.f53d7bf0-3d5b-47a5-bbac-c3f592f6f5ce}
       {v}_{{q + 1}}{ = {v}_{{1}}C_{{q}}{{ \left( {v}_{{1}} \mu  \right) }}^{{q}}\hspace{2mm} , q = 0 , 1 , ...}
       \end{equation}

        \noindent where the $ C_{{r}} ' s $ are the Catalan numbers

       \begin{equation}\label{EQUATION.99d57911-7fa5-4261-8731-8d1e00405620}
       C_{{r}} = \frac {{{ \left( 2r \right) } ! }}{{{ \left( r + 1 \right) } ! r ! }} , r = 0 , 1 , 2 , ... ; \hspace{2mm}\hspace{2mm}C_{{r}} \sim \frac {{4^{{r}}}}{{\sqrt{{ \pi }}{r}^{{3 / 2}}}}\hspace{2mm}for\hspace{2mm}r \rightarrow  \infty .
       \end{equation}

        \noindent For a $  \mu  -  $convolution of order $ m $ with $  \mu  > 0 $ consisting of $ m\hspace{2mm} $initial nonnegative terms $ u_{{1}} , u_{{2}} , ... , u_{{m}}\hspace{2mm} $we define

       \begin{equation}\label{EQUATION.40d9a854-70c1-48f9-8b40-cbb1b5c8e925}
       v_{{1}}\overset{{def.}}{{ = }}{\underset{{k = 1 , 2 , ... , m}}{{max}}{{ \left( \frac {{u_{{k}}}}{{C_{{k - 1}}{ \mu }^{{k - 1}}}} \right) }}^{{1 / k}}.}
       \end{equation}

        \noindent The sequence $ u_{{k}}\hspace{2mm} $is then bounded by the $  \mu  $-convolution $ { \left\{ v_{{k}} \right\} } $ of order 1 and initial term $ v_{{1}} $, i.e. for any $ q :  $

       \begin{equation}\label{EQUATION.e8a58010-8ad4-4865-a3b0-fb663d4c716d}
       u_{{q + 1}} \le v_{{q + 1}} = {v}_{{1}}^{{q + 1}}{{ \mu }}^{{q}}C_{{q}} \sim {v}_{{1}}^{{q + 1}}{ \mu }^{{q}}\frac {{4^{{q}}}}{{\sqrt{{ \pi }}q^{{3 / 2}}}}\hspace{2mm}for\hspace{2mm}q \rightarrow  \infty .
       \end{equation}

        \noindent $ {} $

       \EmptyParagraph
     
    \end{proposition}\begin{proof}
     
        \noindent One of the defining relations for the Catalan numbers is

       \begin{equation}\label{EQUATION.54e872f5-8500-419b-bae1-f1a0700124e9}
       C_{{s}} = C_{{0}}C_{{s - 1}} + C_{{1}}C_{{s - 2}} + ... + C_{{s - 1}}C_{{0}} , \hspace{2mm}s = 1 , 2 , ...
       \end{equation}

        \noindent Equation  \eqref{EQUATION.f53d7bf0-3d5b-47a5-bbac-c3f592f6f5ce} follows then from a straightforward proof by induction and the asymptotic form of  \eqref{EQUATION.99d57911-7fa5-4261-8731-8d1e00405620} results from Stirling's formula  \cite{BIBITEM.921ffd43-b1a4-4972-97d9-34737451af85}. The result of  \eqref{EQUATION.e8a58010-8ad4-4865-a3b0-fb663d4c716d} is a direct consequence of  \eqref{EQUATION.40d9a854-70c1-48f9-8b40-cbb1b5c8e925} which states that $ u_{{k}} \le v_{{k}}\hspace{2mm} $for $ k = 1 , 2 , ... , m. $

       \EmptyParagraph
     
    \end{proof}
     
        \noindent We will now use these results to assess the radius of convergence of the power series $  \sum {{{ \beta }}_{{m}}t^{{m}}.} $
     
       \begin{theorem}\label{THEOREM.13683af9-2bbe-42eb-88a0-a49b06464df4}
     
        \noindent A lower bound for the radius of convergence of $  \sum {{ \beta }}_{{m}}t^{{m}} $ is

       \begin{equation}\label{EQUATION.09ab1bef-b401-469f-a9ff-0ce413fd6467}
       LBRC{ \left(  \alpha  , { \left| a_{{0}} \right| } \right) }\overset{{def.}}{{ = }}\hspace{2mm}\frac {{1}}{{4}}{{ \left( \underset{{k = 1 , 2 , ... , n - 1}}{{max}}{{ \left( \frac {{S_{{k + 1}}{ \left(  \alpha  , { \left| a_{{0}} \right| } \right) { \left\| M \right\| }}}}{{C_{{k - 1}}{{ \left( { \left| a_{{0}} \right| } \right) }}^{{1 - k / n}}}} \right) }}^{{1 / k}} \right) }}^{{ - 1}}.
       \end{equation}

       \EmptyParagraph
     
    \end{theorem}\begin{proof}
     
        \noindent We first define

       \begin{equation}\label{}
       v_{{1}}{ \left(  \alpha  , { \left| a_{{0}} \right| } \right) }\overset{{def.}}{{ = }}\hspace{2mm}\underset{{k = 1 , 2 , ... , n - 1}}{{max}}{{ \left( \frac {{S_{{k + 1}}{ \left(  \alpha  , { \left| a_{{0}} \right| } \right) }}}{{C_{{k - 1}}{ \mu }^{{k - 1}}}} \right) }}^{{1 / k}} , 
       \end{equation}

        \noindent and bear in mind that the $  \mu  -  $convolution $ S_{{q}}{ \left(  \alpha  , { \left| a_{{0}} \right| } \right) }\hspace{2mm} $of order $ n - 1 $ starts at $ q = 2.\hspace{2mm} $ Therefore $ S_{{2}}{ \left(  \alpha  , { \left| a_{{0}} \right| } \right) }\hspace{2mm} $is $ u_{{1}}\hspace{2mm} $of Proposition  \ref{PROPOSITION.bd64019d-cd55-42a1-81d9-61abb499b4ed}, and more generally $ S_{{m}}{ \left(  \alpha  , { \left| a_{{0}} \right| } \right) }\hspace{2mm} $is $ u_{{m - 1}}.\hspace{2mm} $The results of  \eqref{EQUATION.e89e2e3b-0afd-40bb-8df9-8190e5f926b8} and  \eqref{EQUATION.e8a58010-8ad4-4865-a3b0-fb663d4c716d} then yield for any $ m $

       \begin{equation}\label{}
       w{ \left( m \right) } \le S_{{m}}{ \left(  \alpha  , { \left| a_{{0}} \right| } \right) } \le v_{{1}}{{ \left(  \alpha  , { \left| a_{{0}} \right| } \right) }}^{{m - 1}} \times { \mu }^{{m - 2}} \times C_{{m - 2}}\hspace{2mm}
       \end{equation}

       \begin{equation}\label{}
        \sim v_{{1}}{{ \left(  \alpha  , { \left| a_{{0}} \right| } \right) }}^{{m - 1}}{ \mu }^{{m - 2}}\frac {{4^{{m - 2}}}}{{\sqrt{{ \pi }}{{ \left( m - 2 \right) }}^{{3 / 2}}}}\hspace{2mm}for\hspace{2mm}m \rightarrow  \infty .
       \end{equation}

        \noindent The series$ \hspace{2mm} \sum { \left| {{ \beta }}_{{m}}t^{{m}} \right| } $ converges if the sequence $ {{ \left\{ {{{{ \left| { \beta }_{{m}} \right| }}}^{{1 / m}}}{ \left| t \right| } \right\} }}_{{m = 1 , 2 , ...}} $ can be bounded by a sequence converging to a limit $<$1. Given Eq.  \eqref{EQUATION.3e6e9933-f6ec-48a5-9d58-d9402b658d2e} and the fact that $  \rho  = { \left| a_{{0}} |  , \hspace{2mm} \right. } $we have

       \begin{equation}\label{}
       {{{ \left| { \beta }_{{m}} \right| }}}^{{1 / m}}{ \left| t \right| } = { \rho }^{{1 / n - 1 / m}}{{{ \left| W{{ \left( m \right) }}_{{1}} \right| }}^{{1 / m}}{ \left| t \right| }}
       \end{equation}

       \begin{equation}\label{}
        \le { \rho }^{{1 / n - 1 / m}}{{w{ \left( m \right) }}}^{{1 / m}}{ \left| t \right| } \le { \left| t \right| }{ \rho }^{{1 / n - 1 / m}}{{ \left( v_{{1}}{{ \left(  \alpha  , { \left| a_{{0}} \right| } \right) }}^{{m - 1}}{ \mu }^{{m - 2}}C_{{m - 2}} \right) }}^{{1 / m}}
       \end{equation}

       \begin{equation}\label{}
        \sim \frac {{4{ \left| t \right| } \times { \left\| M \|  \times v_{{1}}{ \left(  \alpha  , { \left| a_{{0}} \right| } \right) } \right. }}}{{{{{ \left( { \left| a_{{0}} \right| } \right) }}^{{1 - 1 / n}}}}}\hspace{2mm} = 4{ \left| t \right| }\underset{{k = 1 , 2 , ... , n - 1}}{{max}}{{ \left( \frac {{S_{{k + 1}}{ \left(  \alpha  , { \left| a_{{0}} \right| } \right) }{ \left\| M \right\| }}}{{C_{{k - 1}}{{ \left( { \left| a_{{0}} \right| } \right) }}^{{1 - k / n}}}} \right) }}^{{1 / k}}for\hspace{2mm}m \rightarrow  \infty \hspace{2mm}\hspace{2mm}
       \end{equation}

        \noindent which yields the desired result$ \hspace{2mm} $of  \eqref{EQUATION.09ab1bef-b401-469f-a9ff-0ce413fd6467}.

       \EmptyParagraph
     
    \end{proof}
     
        \noindent We now know that $  \sum {{ \beta }}_{{m}}t^{{m}} $ has a positive radius of convergence which we call $ RC{ \left( a \right) }\hspace{2mm} $to emphasize its dependence on the vector $ a = { \left( a_{{k}} \right) } $ of coefficients. We do not have an analytical expression for $ RC{ \left( a \right) }\hspace{2mm} $but we do have the lower bound $ LBRC \left(  \alpha  , { \left| a_{{0}} \right| } \right)  $.

       For any $ \hspace{2mm}{ \left| t \right| } < RC{ \left( a \right) } $ we may then define

       \begin{equation}\label{}
       c\overset{{def.}}{{ = }}\frac {{t / RC{ \left( a \right) } + 1}}{{2}} < 1
       \end{equation}

        \noindent in which case there exists $ A > 0\hspace{2mm} $such that

       \begin{equation}\label{EQUATION.176c9d9f-d3da-4e04-a2e6-e1cbf6f148b1}
       { \left| b_{{m}}t^{{m}} |  = { \left| K{ \left( m , {B}_{{m}}^{{1}} \right) }t^{{m}} \right| } \right. }{ =  | { \beta }_{{m}} | {{{ \left| t \right| }}}^{{m}} \le {Ac}^{{m}} , \hspace{2mm}m = 1 , 2 , ...}
       \end{equation}

       \EmptyParagraph

       \EmptyParagraph

       In order to show that when the series $  \sum x{ \left( t , k \right) } =  \sum { \beta }_{{m}}{e}^{{2k \pi m \times i / n}}t^{{m}} $ converge$ \hspace{2mm} $they provide the $ n\hspace{2mm} $roots of Equation  \eqref{EQUATION.138330db-ce67-4d93-a8f0-635652033f0e}, we need to prove a result on the growth of the $ K{ \left( d , {B}_{{q + 1}}^{{p}} \right)  ' s.} $

       \EmptyParagraph
     
    \begin{lemma}\label{LEMMA.b4d272cc-9d30-4ac8-83ba-43bfaba948b7}
     
        \noindent If $ {\hspace{2mm}} $$ { \left| t \right| } < RC{ \left( a \right) } $ then there exists $ A > 0\hspace{2mm} $and $ \hspace{2mm}c\hspace{2mm}{ \left( 0 < c < 1 \right) }\hspace{2mm} $such that

       \begin{equation}\label{EQUATION.85d5a972-1074-4adf-9fa1-7f2caff38e8a}
       { \left| K \left( d , {B}_{{q + 1}}^{{p}} \right)  \right| } \le A^{{p}}{{ \left( q + 1 \right) }}^{{p - 1}}{{ \left( {c / { \left| t \right| }} \right) }}^{{d}} , \hspace{2mm} \forall d \ge p \ge 1 , \hspace{2mm}q \ge 0.
       \end{equation}
       \end{lemma}\begin{proof}
     
        \noindent We prove the result by induction on $ p.\hspace{2mm}\hspace{2mm} $Equation  \eqref{EQUATION.176c9d9f-d3da-4e04-a2e6-e1cbf6f148b1} shows that if $ d \le q + 1 $ then

       \begin{equation}\label{}
        \left| K \left( d , {B}_{{q + 1}}^{{1}} \right)  |  = {{ \left| { \beta }_{{d}} \right| } \le A{{ \left( c / { \left| t \right| } \right) }}^{{d}}} \right. 
       \end{equation}

        \noindent which proves  \eqref{EQUATION.85d5a972-1074-4adf-9fa1-7f2caff38e8a} for $ p = 1 ; \hspace{2mm} $ \eqref{EQUATION.85d5a972-1074-4adf-9fa1-7f2caff38e8a} is trivially true if $ d > q + 1\hspace{2mm} $since then $ K{ \left( d , {B}_{{q + 1}}^{{1}} \right) } = 0.\hspace{2mm} $For $ p = 2 $, Eq.  \eqref{EQUATION.6235dca5-7ffb-4919-8483-b028ecf45f0a} yields

       \begin{equation}\label{EQUATION.6cda6834-3396-4d95-b424-8dd51cdb0b73}
       { \left| K{ \left( d , {B}_{{q + 1}}^{{2}} \right) } \right| } = {{ \left| \underset{{m = 1}}{\overset{{min{ \left( d - 1 , q + 1 \right) }}}{ \sum }}b_{{m}}K{ \left( d - m , {B}_{{q + 1}}^{{1}} \right) } \right| } \le }
       \end{equation}

       \begin{equation}\label{}
       \underset{{m = 1}}{\overset{{q + 1}}{ \sum }}{{ \left| b_{{m}} \right| }{{ \left| K{ \left( d - m , {B}_{{q + 1}}^{{1}} \right) } \right| } \le \underset{{m = 1}}{\overset{{q + 1}}{ \sum }}A{{ \left( c / { \left| t \right| } \right) }}^{{m}}A{{ \left( c / { \left| t \right| } \right) }}^{{d - m}}} = {A}^{{2}}{ \left( q + 1 \right) }{{ \left( c / { \left| t \right| } \right) }}^{{d}}.}
       \end{equation}

        \noindent We now assume that  \eqref{EQUATION.85d5a972-1074-4adf-9fa1-7f2caff38e8a} is true up to order $ p < d $ and calculate $  \left| K{ \left.  \left( d , {B}_{{q + 1}}^{{p + 1}} \right)  \right| }.\hspace{2mm} \right.  $Equation  \eqref{EQUATION.6235dca5-7ffb-4919-8483-b028ecf45f0a} yields

       \begin{equation}\label{}
       { \left| K{ \left( d , {B}_{{q + 1}}^{{p + 1}} \right) } \right| } =  | \underset{{m = 1}}{\overset{{min{ \left( d - p , q + 1 \right) }}}{ \sum }}b_{{m}}K{ \left( d - m , {B}_{{q + 1}}^{{p}} \right) } |  \le \underset{{m = 1}}{\overset{{q + 1}}{ \sum }}{ \left| b_{{m}} \right| }{ \left| K{ \left( d - m , {B}_{{q + 1}}^{{p}} \right) } \right| }
       \end{equation}

       \begin{equation}\label{}
       {{ \le \underset{{m = 1}}{\overset{{q + 1}}{ \sum }}A{{ \left( c / { \left| t \right| } \right) }}^{{m}}A^{{p}}{{ \left( q + 1 \right) }}^{{p - 1}}{{ \left( c / { \left| t \right| } \right) }}^{{d - m}}} = A^{{p + 1}}{{ \left( q + 1 \right) }}^{{p}}{{ \left( c / { \left| t \right| } \right) }}^{{d}}}
       \end{equation}

        \noindent which is the desired result.

       \EmptyParagraph
     
    \end{proof}
     
        \noindent In the next section we bring together previous results and prove that when it converges, the power series $ x{ \left( t , k \right) \hspace{2mm}} $provides the roots of Equation  \eqref{EQUATION.138330db-ce67-4d93-a8f0-635652033f0e}.

      \section{Main result}\label{SECTION.b5ae86bd-c55e-42ee-a9cf-82823007b959}
      \begin{theorem}\label{THEOREM.f530226b-2fd6-4f96-8d4d-5307487b9c36}
     
        \noindent We consider the following polynomial equation, parameterized by $ t > 0\hspace{2mm} $and with $ a_{{0}} =  \rho e^{{i \theta }} \ne 0 :  $

       \begin{equation}\label{EQUATION.a70f8bff-9810-4899-9f2e-c5d747660813}
       x^{{n}} = { \left( a_{{n - 1}}x^{{n - 1}} + a_{{n - 2}}x^{{n - 2}} + ... + a_{{1}}x + a_{{0}} \right) }t^{{n}}.
       \end{equation}

        \noindent $ \hspace{2mm} $We define $  \alpha \overset{{def.}}{{ = }}\underset{{k = 1 , 2 , ... , n - 1}}{{max}}{ \left| a_{{k}} \right| }\hspace{2mm} $and recall the definition of the matrix $ M $ of  \eqref{EQUATION.b70e0749-3e15-4609-9d3f-fa3a5aad21ca}. We then define $ n - 1\hspace{2mm} $numbers $ S_{{k}}{ \left(  \alpha  , { \left| a_{{0}} \right| } \right) }\hspace{2mm} $recursively as follows:

       \begin{equation}\label{}
       S_{{2}}{ \left(  \alpha  , { \left| a_{{0}} \right| } \right) }\overset{{def.}}{{ = }} \alpha  +  \alpha  / { \left\| M \|  ,  \right. }
       \end{equation}

       \begin{equation}\label{}
       S_{{q}}{ \left(  \alpha  , { \left| a_{{0}} \right| } \right) }\overset{{def.}}{{ = }}\frac {{1}}{{{{ \left| a_{{0}} \right| }}}}{ \left( \underset{{p = 1}}{\overset{{q - 2}}{ \sum }}{ \left[ { \left\| M \| S_{{p + 1}} \left(  \alpha  , { \left| a_{{0}} \right| } \right)  +  \alpha  \right. } \right] }S_{{q - p}}{ \left(  \alpha  , { \left| a_{{0}} \right| } \right) } \right)  +  \alpha  +  \alpha  / { \left\| M \right\| }{ , q = 3 , 4 , ... , n}}
       \end{equation}

       \EmptyParagraph

       We also define recursively the infinite sequence of $ { \left( n - 1 \right) }\hspace{2mm} $dimensional vectors $ W{ \left( q \right) } :  $

       \begin{equation}\label{}
       W{ \left( 1 \right) } = a_{{0}}{ \left( \begin{matrix}1 & 1 & ... & ... & 1\end{matrix} \right) } ' \hspace{2mm} ; \hspace{2mm}W{ \left( 2 \right) } = \frac {{a_{{1}}}}{{n}}{ \left( \begin{matrix}1 & 2 & ... & ... & n - 1\end{matrix} \right) } ' 
       \end{equation}

       \begin{equation}\label{}
       W{ \left( q \right) } = \frac {{1}}{{{a_{{0}}}}} \left[  - M\underset{{p = 1}}{\overset{{q - 2}}{ \sum }}{ \left[ u.W{ \left( p + 1 \right) } \right] }W{ \left( q - p \right) } + U\underset{{p = 1}}{\overset{{min{ \left( n - 1 , q - 1 \right) }}}{ \sum }}{ \left[ A_{{p}}.W{ \left( q - p \right) } \right] } \right]  , q = 3 , 4 , ...
       \end{equation}

        \noindent where:

       1. u is the (n-1)-dimensional row vector having 1 in first position and zeros elsewhere.

       2. U is the (n-1)-dimensional column vector $ \frac {{1}}{{n}}{ \left( \begin{matrix}1 & 2 & ... & ... & n - 1\end{matrix} \right)  ' }. $

       3. Each $ A_{{p}}\hspace{2mm} $(p=1,2,...,n-1) is the $ { \left( n - 1 \right) } $-dimensional row vector with $ a_{{p}} $ in

       $ p - th $ position and zeros elsewhere.

       If we define the sequence

       \begin{equation}\label{EQUATION.ec208c74-9ac4-4db2-8834-278e169b0271}
       { \beta }_{{m}} = { \rho }^{{m / n}} \times e^{{i{ \left( m \theta  \right)  / n}}} \times W{{ \left( m \right) }}_{{1}} / a_{{0}} , \hspace{2mm}m = 1 , 2 , ...
       \end{equation}

        \noindent then $  \sum { \left| { \beta }_{{m}}t^{{m}} | \hspace{2mm} \right. } $converges for $ { \left| t |  < RC{ \left( a \right) }\hspace{2mm} \right. } $and in particular for

       \begin{equation}\label{EQUATION.1dc121b3-89ac-4ca0-a412-ddf4d732928c}
       {{ \left| t \right| } < LBRC{ \left(  \alpha  , { \left| a_{{0}} \right| } \right) }\overset{{def.}}{{ = }}\hspace{2mm}\frac {{1}}{{4}}{{ \left( \underset{{k = 1 , 2 , ... , n - 1}}{{max}}{{ \left( \frac {{S_{{k + 1}}{ \left(  \alpha  , { \left| a_{{0}} \right| } \right) { \left\| M \right\| }}}}{{C_{{k - 1}}{{ \left( { \left| a_{{0}} \right| } \right) }}^{{1 - k / n}}}} \right) }}^{{1 / k}} \right) }}^{{ - 1}}.}
       \end{equation}

        \noindent If $|$$ t |  < RC{ \left( a \right) } $ the $ n $ convergent series

       \begin{equation}\label{EQUATION.6bba582d-2e42-4b60-8a07-25ad269ae181}
       x{ \left( t , k \right) } = \underset{{m = 1}}{\overset{{ \infty }}{ \sum }}{ \beta }_{{m}}e^{{2k \pi m \times i / n}} \times t^{{m}} , \hspace{2mm}k = 0 , 1 , ... , n - 1
       \end{equation}

        \noindent are the n roots of Eq.  \eqref{EQUATION.a70f8bff-9810-4899-9f2e-c5d747660813} in the sense that for each $ k = 0 , 1 , ... , n - 1 $ the expression

       \begin{equation}\label{EQUATION.c46a0cad-9f5d-4cfc-9b0b-a3f16336a2e4}
       E{ \left( x{{ \left( t , k \right) }}_{{q}} \right) }\overset{{def.}}{{ = }}{{ \left[ x{{ \left( t , k \right) }}_{{q}} \right] }}^{{n}} - { \left( {a}_{{{n - 1}}}{{ \left[ x{{ \left( t , k \right) }}_{{q}} \right] }}^{{n - 1}} + a_{{n - 2}}{{ \left[ x{{ \left( t , k \right) }}_{{q}} \right] }}^{{n - 2}} + ... + a_{{1}}x{{ \left( t , k \right) }}_{{q}} + {a}_{{0}} \right) }t^{{n}}
       \end{equation}

        \noindent approaches $ 0 $ for $ q \rightarrow  \infty  $. (The $ x{{ \left( t , k \right) }}_{{q}} ' s $ are the partial sums up to order $ q\hspace{2mm} $(Eq.  \eqref{EQUATION.6e641714-973c-4278-a36b-890fb646588e} ). 
     
       \end{theorem}\begin{proof}
     
        \noindent Given Eq.  \eqref{EQUATION.26e6c154-9a63-44c9-ae45-15c4a1a32ad6} the polynomial in $ t\hspace{2mm} $on the right-hand side of  \eqref{EQUATION.c46a0cad-9f5d-4cfc-9b0b-a3f16336a2e4} only has terms of the form $ t^{{n + q + s}}\hspace{2mm} $with $ s \ge 0.\hspace{2mm} $We will first show that $ {{ \left[ x{{ \left( t , k \right) }}_{{q}} \right] }}^{{n}} $ contributes a series of such terms that approaches $ 0\hspace{2mm} $when $ q \rightarrow  \infty .\hspace{2mm} $Indeed,  \eqref{EQUATION.85d5a972-1074-4adf-9fa1-7f2caff38e8a} shows that the modulus of the contribution of $ {{ \left[ x{{ \left( t , k \right) }}_{{q}} \right] }}^{{n}}\hspace{2mm} $to the terms $ t^{{n + q + s}}\hspace{2mm} $is

       \begin{equation}\label{}
       {{{ \left| \underset{{s = 0}}{\overset{{q{ \left( n - 1 \right) } - n}}{ \sum }}t^{{n + q + s}}K{ \left( n + q + s , {B}_{{q}}^{{n}} \right) } \right| } \le \underset{{s = 0}}{\overset{{ \infty }}{ \sum }}{{ \left| t \right| }}^{{n + q + s}}A^{{n}}q^{{n - 1}}{{ \left( c / { \left| t \right| } \right) }}^{{n + q + s}}} = }
       \end{equation}

       \begin{equation}\label{}
       A^{{n}}{q}^{{n - 1}}c^{{n + q}}\underset{{s = 0}}{\overset{{ \infty }}{ \sum }}{c}^{{s}}
       \end{equation}

        \noindent which approaches $ 0 $ when $ q \rightarrow  \infty . $

       We similarly consider the modulus of the contribution of each $ {{ \left[ x{{ \left( t , k \right) }}_{{q}} \right] }}^{{p}}\hspace{2mm} $ of  \eqref{EQUATION.c46a0cad-9f5d-4cfc-9b0b-a3f16336a2e4} (for $  \left. 1 \le p \le n - 1 \right) \hspace{2mm} $to the series of terms $ t^{{n + q + s}}\hspace{2mm}.\hspace{2mm} $Given the term $ t^{{n}}\hspace{2mm} $appearing on the right-hand side of  \eqref{EQUATION.c46a0cad-9f5d-4cfc-9b0b-a3f16336a2e4}, this modulus is

       \begin{equation}\label{}
       { \left| \underset{{s = 0}}{\overset{{q{ \left( p - 1 \right) }}}{ \sum }}t^{{q + s}}K{ \left( q + s , {B}_{{q}}^{{p}} \right) } \right| } \le \underset{{s = 0}}{\overset{{ \infty }}{ \sum }}{{ \left| t \right| }}^{{q + s}}A^{{p}}q^{{p - 1}}{{ \left( c / { \left| t \right| } \right) }}^{{q + s}} = 
       \end{equation}

       \begin{equation}\label{}
       A^{{p}}{q}^{{p - 1}}c^{{q}}\underset{{s = 0}}{\overset{{ \infty }}{ \sum }}{c}^{{s}}
       \end{equation}

        \noindent which also approaches $ 0 $ when $ q \rightarrow  \infty .\hspace{2mm} $(Because $ a_{{0}} = {b}_{{1}}^{{n}} $ the term $ a_{{0}}t^{{n}}\hspace{2mm} $on the right-hand side of  \eqref{EQUATION.c46a0cad-9f5d-4cfc-9b0b-a3f16336a2e4} cancels out with the term $ {b}_{{1}}^{{n}}t^{{n}}\hspace{2mm} $appearing in the expansion of $ {{ \left[ x{{ \left( t , k \right) }}_{{q}} \right] }}^{{n}}. $) This completes the proof since each one of the $ {{ \left[ x{{ \left( t , k \right) }}_{{q}} \right] }}^{{p}}\hspace{2mm}\hspace{2mm} $terms on the right-hand side of $ E{ \left( x{{ \left( t , k \right) }}_{{q}} \right) }\hspace{2mm} $approaches $ 0 $ when $ q \rightarrow  \infty . $

       \EmptyParagraph
     
    \end{proof}
     
        \noindent $ \hspace{2mm} $

       {\bf Remarks:}

       1. The functions $ x{ \left( t , k \right) } $ of Eq.  \eqref{EQUATION.6bba582d-2e42-4b60-8a07-25ad269ae181} can be viewed as Taylor expansions in the variable $ t\hspace{2mm} $of the roots of Equation  \eqref{EQUATION.a70f8bff-9810-4899-9f2e-c5d747660813}: when $ t $ is small the first few terms of the series provide approximate values for the roots (a numerical example is given below).

       2. In view of  \eqref{EQUATION.ec208c74-9ac4-4db2-8834-278e169b0271} the radius of convergence $ RC{ \left( a \right) }\hspace{2mm} $of $ x{ \left( t , u \right) } $ as a power series in $ t\hspace{2mm} $can be calculated numerically as

       \begin{equation}\label{}
       RC{ \left( a \right) } = \underset{{m \rightarrow  \infty }}{{lim\hspace{2mm}inf\hspace{2mm}}}Q{ \left( m \right) }
       \end{equation}

        \noindent where

       \begin{equation}\label{}
       Q{ \left( m \right) }\overset{{def.}}{{ = }}\frac {{1}}{{{{ \left( { \left| a_{{0}} \right| } \right) }}^{{1 / n}}{{{ \left( { \left| W{{ \left( m \right) }}_{{1}} \right| } \right) }}}_{{}}^{{1 / m}}}}
       \end{equation}

       \EmptyParagraph

       3. Each term $ S_{{q}}{ \left(  \alpha  , { \left| a_{{0}} \right| } \right) } $ is a polynomial in $  \alpha \hspace{2mm} $with no constant term and positive coefficients. This insures that $ LBRC{ \left(  \alpha  , { \left| a_{{0}} \right| } \right) }\hspace{2mm}\hspace{2mm} $tends to infinity for $  \alpha  \rightarrow 0 $ . Therefore  \eqref{EQUATION.1dc121b3-89ac-4ca0-a412-ddf4d732928c} is satisfied for $  \alpha \hspace{2mm} $small enough which shows that the functions $ x{ \left( t , u \right) }\hspace{2mm} $provide the roots when the coefficients $ a_{{k}}\hspace{2mm} $other than $ a_{{0}} $ are small enough.

       4. Each $ S_{{q}}{ \left(  \alpha  , { \left| a_{{0}} \right| } \right) } $ is a decreasing function of $ { \left| a_{{0}} | \hspace{2mm} \right. } $that approaches $  \alpha  +  \alpha  / { \left\| M \right\| } $ when $ { \left| a_{{0}} \right| } \rightarrow  \infty . $ With $ {{ \left( { \left| a_{{0}} \right| } \right) }}^{{1 - k / n}} $ in the denominator of  \eqref{EQUATION.1dc121b3-89ac-4ca0-a412-ddf4d732928c}, the function $ LBRC{ \left(  \alpha  , { \left| a_{{0}} \right| } \right) }\hspace{2mm} $ tends to infinity for $ { \left| a_{{0}} |  \rightarrow  \infty .\hspace{2mm} \right. } $Therefore the functions $ x{ \left( t , u \right) }\hspace{2mm} $provide the roots when $ {{ \left| a \right. }}_{{0}} | \hspace{2mm} $is large enough.

       \EmptyParagraph

       As a numerical illustration we consider the equation

       \begin{equation}\label{EQUATION.902d8e88-2960-4aec-b534-fd3c13e6f2c3}
       x^{{6}} = { \left(  - x^{{5}} + x^{{4}} - 2x^{{3}} - 3x^{{2}} + 2x + 8 \right) }t^{{6}}
       \end{equation}

        \noindent with $ t $ set equal to 1. The lower bound $ LBRC{ \left(  \alpha  , { \left| a_{{0}} \right| } \right) \hspace{2mm}\hspace{2mm}} $for the radius of convergence of the power series $ x{ \left( t , u \right) }\hspace{2mm} $ is 0.094. Therefore we do not know whether $ x{ \left( t , u \right) } $ converges with $ t = 1.\hspace{2mm} $$ {} $However $ LBRC{ \left(  \alpha  , { \left| a_{{0}} \right| } \right) \hspace{2mm}} $is an extremely conservative bound calculated only with $ { \alpha  , { \left| a_{{0}} \right| }\hspace{2mm}} $and $ { \left\| M \right\| } $; $ LBRC{ \left(  \alpha  , { \left| a_{{0}} \right| } \right) \hspace{2mm}} $thus represents a worst-case scenario on the coefficients $ a_{{k}} $. For these reasons the condition $ { \left| t \right| } < LBRC{ \left(  \alpha  , { \left| a_{{0}} \right| } \right) }{\hspace{2mm}} $ of  \eqref{EQUATION.1dc121b3-89ac-4ca0-a412-ddf4d732928c} is probably of limited practical use. Its usefulness lies in the qualitative Remarks 3 and 4 above, which state that the function $ x{ \left( t , u \right) } $ does provide the roots for $  \alpha \hspace{2mm} $small enough or $ { \left| a_{{0}} | \hspace{2mm} \right. } $large enough. 
     
       \begin{figure}[H]\centering \includegraphicsX{scale=1}{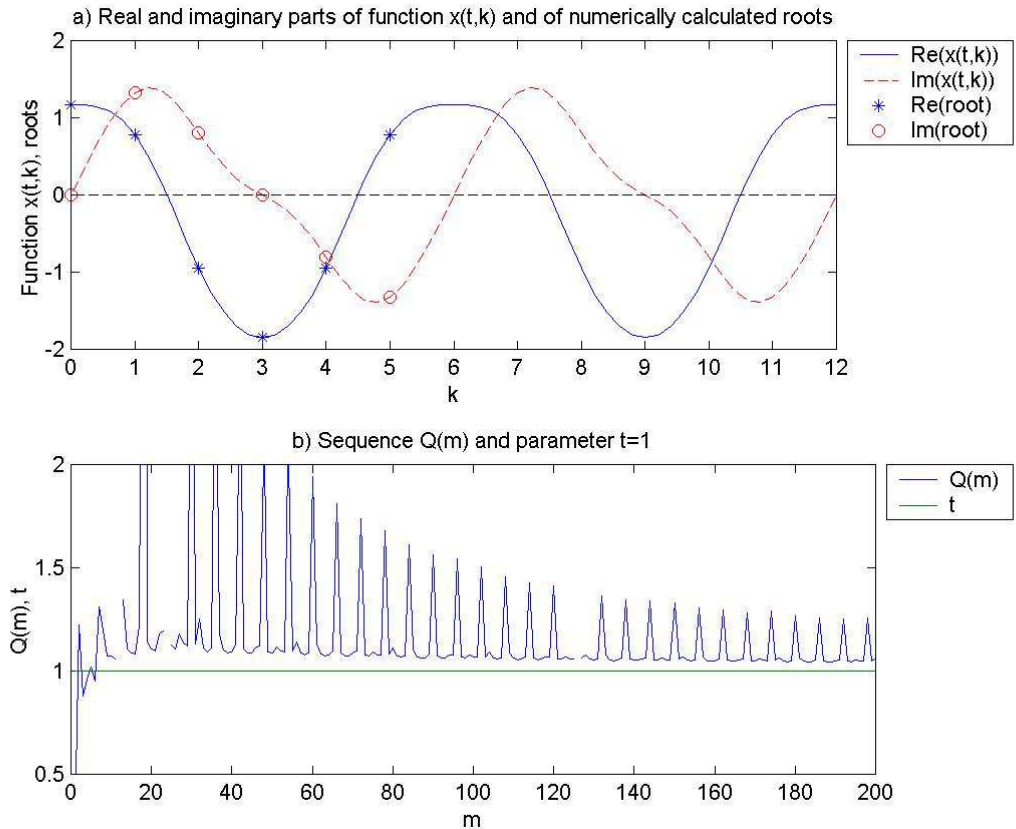}\caption{Function x(t,k) (with t=1) and sequence Q(m)}
      \end{figure}

       \EmptyParagraph

       From Figure 5.1b the radius of convergence $ RC{ \left( a \right) } $ assessed numerically as the limit of $ Q{ \left( m \right) } $ appears to be about 1.05. Therefore the power series $ x{ \left( t , u \right) } $ does converge with $ t = 1. $ (A few missing values in the $ Q{ \left( m \right) }\hspace{2mm} $sequence arise when $ W{{ \left( m \right) }}_{{1}} = 0 $, i.e. $ Q{ \left( m \right) } =  \infty . $)

       The real and imaginary parts of $ x{ \left( t , k \right) }\hspace{2mm} $are plotted in Figure 5.1a over two periods of length $ n = 6 $ (200 terms are used to calculate the series). These two parts taken at $ t = 1\hspace{2mm} $and $ k = 0 , 1 , 2 , 3 , 4 , 5\hspace{2mm} $ coincide with those obtained using Matlab's built-in polynomial equation routine to solve Eq.  \eqref{EQUATION.902d8e88-2960-4aec-b534-fd3c13e6f2c3} with $ t = 1\hspace{2mm} $(stars and circles in Figure 5.1a).

       With real coefficients for the polynomial equation the complex roots come in conjugate pairs. There are two such pairs. One for $ k = 1 , k = 5 $ and the other for $ k = 2 , k = 4 $. In addition there are two real roots at $ k = 0\hspace{2mm} $and at $ k = 3.\hspace{2mm} $In this and some other cases these real roots appear to be at local minima or maxima of the $ Re{ \left( x{ \left( t , k \right) } \right) }\hspace{2mm} $function. These special behaviors of and relationships between the real and imaginary parts no doubt arise from particular patterns in the sequence $ { \left\{ W{{ \left( q \right) }}_{{1}} \right\} } $ which have yet to be explored. These structures disappear with complex coefficients since in this case there are no more complex conjugate roots.

       The first five terms of the series $ x{ \left( t , u \right) }\hspace{2mm} $are

       \[
       x{ \left( t , u \right) } \approx \sqrt{{2}}e^{{2u \pi  \times i / 6}}t + 0.08333e^{{4u \pi  \times i / 6}}t^{{2}} - 0.18414e^{{6u \pi  \times i / 6}}t^{{3}}
       \]

       \begin{equation}\label{EQUATION.fe2e69ed-e243-4e37-ad14-6318fe0edfb3}
        - 0.14506e^{{8u \pi  \times i / 6}}\hspace{2mm}t^{{4}} + 0.11441e^{{10u \pi  \times i / 6}}t^{{5}}.
       \end{equation}

        \noindent There are not enough terms to calculate the roots when $ t = 1. $ With $ t = 0.4 $ however, the values provided by  \eqref{EQUATION.fe2e69ed-e243-4e37-ad14-6318fe0edfb3} are very close to those calculated numerically (Table  \ref{TABLE.4938efc6-b082-4572-9fb7-e5725be5726f}).

       \EmptyParagraph

                 \begin{table}[H]
                 \begin{center}
                  \caption{x(t,k) and numerically calculated roots of Equation (5.16) with t=0.4}\label{TABLE.4938efc6-b082-4572-9fb7-e5725be5726f}
          {\renewcommand{\arraystretch}{1.1} 
          \tableX{|l|l|l|l|l|l|l|}{|X|X|X|X|X|X|X|}{\hline
        &  \multicolumn{1}{|c|}{k=0} &  \multicolumn{1}{|c|}{k=1} &  \multicolumn{1}{|c|}{k=2} &  \multicolumn{1}{|c|}{k=3} &  \multicolumn{1}{|c|}{k=4} &  \multicolumn{1}{|c|}{k=5 } \\ \hline
      x(0.4,k) &  {\small 0.565} &  {\small 0.290+0.504i} &  \multicolumn{1}{|c|}{{\small -0.300+0.474i}} &  \multicolumn{1}{|c|}{{\small -0.545}} &  \multicolumn{1}{|c|}{{\small -0.300-0.474i}} &  \multicolumn{1}{|c|}{{\small 0.290-0.504i}} \\ \hline
      num. root &  {\small 0.564} &  {\small 0.290+0.504i} &  {\small -0.301+0.474i} &  {\small -0.546} &  {\small -0.301-0.474i} &  {\small 0.290-0.504i} \\ \hline
      }}  
        \end{center}
      \end{table}

       \EmptyParagraph

       If the constant term $ a_{{0}} = 8\hspace{2mm} $in Eq.  \eqref{EQUATION.902d8e88-2960-4aec-b534-fd3c13e6f2c3} is less than approximately 4, the radius of convergence for the series $ x{ \left( t , u \right) }\hspace{2mm} $drops below 1. The method can no longer be used to solve Eq.  \eqref{EQUATION.902d8e88-2960-4aec-b534-fd3c13e6f2c3} with $ t = 1 $. Several approaches have been tried in order to extend the method to the case $ a_{{0}} $ small in which at least one root becomes close to 0. One possibility would be to transform the unknown in such a way that small roots are moved away from 0. For example one could write the polynomial equation in terms of a changed unknown $ y = 1 / x $: if a root $ x $ is small then the corresponding $ y\hspace{2mm} $is large. However this approach did not change the problem. Another more promising possibility would be to inject the parameter $ t $ differently into the equation. For example one could use a similar approach after multiplying the left side of Eq.  \eqref{EQUATION.07dd3654-31ac-4ca3-a570-c25959e51126} by $ t^{{n}}\hspace{2mm} $instead of the right side. Or one could multiply each term of the equation by a power $ t^{{p}}\hspace{2mm} $with a different and well-chosen $ p $ for each term. To date such attempts have proved largely inconclusive.

       Ours is only a first step, which shows that the roots of a particular class of polynomial equations can be expressed explicitly with an infinite number of rational operations and root extractions. It is to be hoped that some variant of the family of periodic functions $ x{ \left( t , u \right) }\hspace{2mm} $ will eventually emerge to provide closed-form expressions for the roots of arbitrary polynomial equations.

\clearpage 
\pagebreak

       \EmptyParagraph

\begin{center}
        \noindent REFERENCES
\end{center}

       \EmptyParagraph


\begin{thebibliography}{XX}
   \bibitem{BIBITEM.32b8409e-4efa-45ca-ae05-adb468ace982}
     
        \noindent N.H. Abel. "Beweis der Unmöglichkeit, algebraische Gleichungen von höheren Graden als dem vierten allgemein aufzulösen." J. reine angew. Math. 1, 65, 1826. Reprinted in Abel, N. H. Oeuvres Completes (Ed. L. Sylow and S. Lie). New York: Johnson Reprint Corp., pp. 66-87, 1988.
     
       \bibitem{BIBITEM.8d478637-3e3c-4e3c-b709-d4353976b5ee}
     
        \noindent G. Belardinelli, Fonctions hypergéométriques de plusieurs variables et résolution analytique des équations algébrique générales.\emph{ Mémoire des Sci. Math.} 145, 1960.
     
       \bibitem{BIBITEM.09654ad0-cb4d-4f5e-97a6-455e481aaa20}
     
        \noindent R. Birkeland. Résolution de l'equation algébrique générale par des fonctions hypergéométriques supérieures, Comp. Rend. 171, 778 ,1921. 
     
       \bibitem{BIBITEM.1e4a0b01-ced4-402b-a8fa-48c735e80cb1}
     
        \noindent Eisenstein, G. Allgemeine Auflösung der Gleichungen von den ersten vier Graden, \emph{J. reine angew. Math.} 27, 81, 1844. 
     
       \bibitem{BIBITEM.2ce5749f-bb1c-4c48-9c60-73b87c19a339}
     
        \noindent P. L. Chebychev, Calculations of the roots of an equation (in Russian), in P. L. Chebychev: Collected Works, v. V, Akademii Nauk SSSR, Moscow, 1951.
     
       \bibitem{BIBITEM.54e711a8-276d-4617-bd00-cf34c8750ced}
     
        \noindent P. A. Lambert, On the solution of algebraic equations in infinite series, \emph{Bull. Am. Math. Soc}. XIV, 467, 1908. 
     
       \bibitem{BIBITEM.921ffd43-b1a4-4972-97d9-34737451af85}
     
        \noindent Weisstein E.W. \emph{et al.} "Catalan Number." From MathWorld--A Wolfram Web Resource. http://mathworld.wolfram.com/CatalanNumber.html. 
     
       
    
       \EmptyParagraph
     
    
    
       \EmptyParagraph
     
    
     \end{thebibliography}
\end{document}